\documentclass[11pt]{amsart}

\usepackage{amsmath,amssymb}
\usepackage{mathrsfs} 
\usepackage{enumerate}
\usepackage{appendix}
\usepackage{wrapfig}

\usepackage{tablefootnote}

\usepackage{pgf,tikz}
\usetikzlibrary{arrows}
\usepackage[all]{xy}
\usepackage{appendix}

\usepackage{float}

\usepackage{array,booktabs,arydshln,xcolor}

\usetikzlibrary {positioning}
\usetikzlibrary{arrows}
\usetikzlibrary{positioning}
\usepackage{multirow}
\usepackage{tikz-cd}

    \usetikzlibrary{patterns.meta}
\tikzstyle{stuff_fill}=[rectangle,draw=black,fill=white,minimum size=1em]
\tikzstyle{stuff_fill_blank}=[rectangle,fill=white,minimum size=1em]

\usetikzlibrary{decorations.pathreplacing}

\usepackage{multicol}

\usepackage{geometry} 
\usepackage{tikz-cd}
\tikzcdset{arrow style=tikz,diagrams={>=stealth}}

\RequirePackage{hyperref}
\hypersetup{pdfpagemode={UseOutlines},
bookmarksopen=true,
bookmarksopenlevel=0,
hypertexnames=false,
colorlinks=true, 
citecolor=teal, 
linkcolor=blue, 
urlcolor=magenta, 
pdfstartview={FitV},
unicode,
breaklinks=true,
}

\usepackage{array,multirow}

\newtheorem{theorem}{Theorem}[section]
\newtheorem{theoremintro}{Theorem}
\newtheorem{corollaryintro}[theoremintro]{Corollary}

\newtheorem{lemma}[theorem]{Lemma}
\newtheorem{proposition}[theorem]{Proposition}
\newtheorem{corollary}[theorem]{Corollary}

\newtheorem*{theorem*}{Theorem}
\newtheorem*{ques*}{Question}
\newtheorem*{prop*}{Proposition}

\theoremstyle{definition}
\newtheorem{definition}[theorem]{Definition}
\newtheorem{example}[theorem]{Example}

\newtheorem{ques}[theorem]{Question}

\newtheorem*{definition*}{Definition}

\theoremstyle{remark}
\newtheorem{remark}[theorem]{Remark}

\numberwithin{equation}{section}

\title[Sublinear bilipschitz equivalence, quasiisometry and Lie groups]{Sublinear bilipschitz equivalence and the quasiisometric classification of solvable Lie groups}

\author{Ido Grayevsky}
\address{Department of Mathematics, Ben Gurion University of the Negev
P.O.B. 653, Be'er Sheva 8410501, Israel}
\email{idogra@post.bgu.ac.il}

\author{Gabriel Pallier}
\address{Univ. Lille, CNRS, UMR 8524 - Laboratoire Paul Painlevé, F-59000 Lille, France}
\email{gabriel.pallier@univ-lille.fr}

\thanks{The first author was partially funded by Israel Science Foundation grant ISF 2919/19 of Yair Glasner, and grant ISF 2990/21
of Ori Parzanchevski. The second author was funded by the DFG 281869850 (RTG 2229). We gratefully acknowledge the financial support from the ConYS programme of the Karlsruhe House for Young Scientists at the KIT, and from the PEPS JCJC programme of the INSMI}

\thanks{}


\date{\today}
\dedicatory{}

\begin{document}

\begin{abstract}
    We prove a product theorem for sublinear bilipschitz equivalences which generalizes the classical work of Kapovich, Kleiner and Leeb on quasiisometries between product spaces. We employ our product theorem to distinguish up to quasiisometry certain families of solvable groups which share the same dimension, cone-dimension and Dehn function; actually we do this by distinguishing them up to sublinear bilipschitz equivalence, which is slightly stronger. As an application, we recover the fact, recently obtained by Bourdon and Rémy with different groups, that there exists uncountably many quasiisometry classes of indecomposable, non-unimodular, high rank solvable Lie groups.
\end{abstract}
\maketitle

\tableofcontents

\section{Introduction}

This study is motivated by the quasiisometric classification of connected Lie groups.
Such a classification would be complete if one could classify the completely solvable Lie groups up to quasiisometry, as every connected Lie group $G$ is quasiisometric to a completely solvable\footnote{The completely solvable Lie groups are the closed subgroups of the upper triangular real matrix group; they are also called real-triangulable, or split-solvable.} group $\rho_0(G)$, called the trigshadow of $G$ \cite{CornulierDimCone}.

Cornulier conjectured that two quasiisometric completely solvable Lie groups should be isomorphic {\cite[Conjecture 19.113]{CornulierQIHLC}}. This is currently open even within the smaller class of simply connected nilpotent groups.

The process of going from $G$ to $\rho_0(G)$ is a reduction procedure. 
Cornulier went further in the reduction procedure and defined two subclasses of completely solvable Lie groups, $(\mathcal C_1)$ and $(\mathcal C_\infty)$, such that every connected Lie group $G$ is $O(\log)$-bilipschitz equivalent to some group $\rho_1(G)$ in $(\mathcal C_1)$, and $O(u)$-bilipschitz equivalent to some group $\rho_\infty(G)$ in $(\mathcal C_\infty)$ for some explicit sublinear function $u$ depending on $G$ (one has $(\mathcal C_\infty) \subset (\mathcal C_1)$, $\rho_\infty \circ \rho_1 = \rho_\infty$ and $\rho_1 \circ \rho_0 = \rho_1$) \cite{CornulierCones11}.
We will recall the definition of $O(u)$-bilipschitz equivalence and Cornulier's reductions farther; let us only specify here that in this language, quasiisometry is $O(1)$-bilipschitz equivalence, that $O(\log)$ and $O(u)$-bilipschitz equivalence are weaker than quasiisometry, and that when $G$ is nilpotent $\rho_\infty(G)$ is the graded nilpotent group associated to the lower central filtration of $G$, which is known to be a quasiisometry invariant by the work of Pansu.

\begin{theorem}[\cite{PanCBN, PansuCCqi}]
    Let $G$ and $G'$ be quasiisometric simply connected nilpotent Lie groups. Then $\rho_\infty(G)$ and $\rho_\infty(G')$ are isomorphic.
\end{theorem}

Our purpose here is to demonstrate that the reduction procedure, in the present case at the level of $\rho_1$, has applications to the quasiisometric classification in the class of Lie groups of exponential growth, which is disjoint from that of nilpotent groups. This stems from the following fact, relying on \cite{CornulierCones11}: since the group $\rho_1(G)$ in the class $(\mathcal C_1)$ have less complicated structure than $G$, it can happen that it splits in a direct product while $G$ did not.

\subsection{Main results}
Given a map between product spaces $f:X_1 \times X_2 \rightarrow Y_1 \times Y_2$, when does it descend to a map between the factors? In~\cite{KKL}, Kapovich, Kleiner and Leeb provide topological conditions on the factors to ensure that this happens for $f$  a homeomorphism. They then deduce similar statements for quasi-isometries, i.e.\ $O(1)$-bilipschitz equivalences,  assuring that quasiisometries preserve the product structure up to a bounded error. Our main contribution is to show that, in a similar setting, $O(u)$-bilipschitz equivalence also descend to the factors and that the product structure is preserved up to a sublinear error term.

For a pointed metric space $(X,x_0)$ we denote by  $|x|=d(x,x_0)$ the \emph{size of $x\in X$}. A \emph{sublinear function} is any real function $u\colon \mathbf R_{\geqslant 0} \to \mathbf R_{\geqslant 0}$ such that $\lim_{r\to\infty}\frac{u(r)}{r}=0$.
Two maps $f$ and $g$ between pointed metric spaces are $O(u)$-close if there is $u_0\in O(u)$ such that $d\big(f(x), g(x)\big)\leq u_0(|x|)$, and \emph{sublinearly close} if there exists such a $u$ that is sublinear.

\begin{definition}[Commuting up to sublinear error]
\label{def:commuting-sublinear-error}
Let $n$ be a positive integer.
    Let $(X_i)$ and $(Y_i)$ be families of pointed metric spaces, $i=1,\ldots, n$.
    Given $i \in \{1, \ldots, n \}$ we say that the diagram
    \begin{equation}
    \label{eq:commutes-up-to-v}
        \begin{tikzcd}
\prod_{i=1}^n X_i \ar[r, "\phi"] \ar[d,"\pi_i"] & \prod_{i=1}^n Y_i \ar[d,"\pi_i"] \\
X_i \ar["\phi_i",r] & Y_i
\end{tikzcd}
    \end{equation}
commutes up to sublinear error if there exists a sublinear function $v$ such that $\phi_i \circ \pi_i$ and $\pi_i \circ \phi$ are $O(v)$-close. For  $x=(x_1,\dots,x_n)$, we denote the product map
$$\Phi(x)= \left( \prod_{i=1}^n\phi_i \right) (x):= \big(\phi_1(x_1),\dots,\phi_n(x_n)\big).$$ 
When the diagrams \eqref{eq:commutes-up-to-v} commute up to error $v$ for all $i$, $\phi$ and $\Phi$ are $O(v)$-close.
\end{definition}

We refer to Section~\ref{sec: product theorem} for the definitions of certain terms in the statement below. Key examples of metric spaces of coarse type I are simply connected Riemannian manifolds with sectional curvature bounded above by a negative constant; examples of spaces of coarse type II are irreducible Riemannian symmetric spaces of noncompact type and high rank.

\begin{theoremintro}\label{thm:KKLgeneralized}
Let $X = M \times \prod_{i=1}^n X_i$ and $Y = N \times \prod_{i=1}^m Y_j$, where $X_i$ and $Y_j$ are of coarse type I or II in the sense of Kapovich, Kleiner and Leeb \cite{KKL}, and $M$ and $N$ are geodesic spaces with asymptotic cones homeomorphic to $\mathbf R^p$ and $\mathbf R^q$ respectively. Let $u$ be a subadditive nondecreasing sublinear function, and let $\phi\colon X \to Y$ be an $O(u)$-bilipschitz equivalence. 
Then $p=q$, $n=m$, and there is an $O(u)$-bilipschitz equivalence $\overline{\phi}:\overline{X}=\prod X_i\rightarrow \prod Y_i=\overline{Y}$ such that the maps $\overline{\pi}\circ\phi$ and $\overline{\phi}\circ\overline{\pi}$ are $O(v)$-close for some subadditive nondecreasing sublinear function $\overline v$. 

Moreover, there exists a bijection $\sigma: \{ 1, \ldots, n \} \to \{ 1, \ldots, n \}$ and, for every $i \in \{1,\ldots, n\}$, an $O(u)$-bilipschitz equivalence $\phi_i: X_i \to Y_{\sigma(i)}$ such that the diagram
\[ \begin{tikzcd}
M \times \prod_{i=1}^n X_i \ar[r, "\overline{\phi}"] \ar[d,"\pi_i"] &  N \times \prod_{i=1}^n Y_i \ar[d,"\pi_{\sigma(i)}"] \\
X_i \ar["\phi_i",r] & Y_{\sigma (i)}
\end{tikzcd} \]
commutes up to a sublinear error $v'$ that is subadditive and nondecreasing. 
\end{theoremintro}

The case $u=1$ follows from \cite[Theorem B]{KKL}; in this case $v$ and $v'$ can be shown to be bounded (see Remark~\ref{remark: O(1) vs O(u) error term}); a more elementary proof in the case $n=1$ and $X_1$ and $Y_1$ are non-elementary regular trees was given in \cite{SoucheWiest}.

We state a first application below, which uses some of the results of \cite{PalSBErankone, GrayevskyRigidity}. By pluriisometry between symmetric spaces we mean that we allow a rescaling in each factor of the de Rham decomposition.

\begin{corollaryintro}
\label{cor:KKL-symmetric}
Let $X$ and $Y$ be two Riemannian globally symmetric spaces with no compact factors. If $X$ and $Y$ are sublinear bilipschitz equivalent, then $X$ and $Y$ are pluriisometric.
\end{corollaryintro}

Corollary~\ref{cor:KKL-symmetric} restricted to quasiisometries instead of sublinear bilispchitz equivalences is due to Kleiner-Leeb \cite[Corollary 1.1.4]{KleinerLeebQI} and builds on the work of Mostow \cite{Mostow70} for the part concerning the rank one factors.

The proof of Theorem~\ref{thm:KKLgeneralized} relies on Theorem~\ref{thm:KKLgeneralizedNEW}, which is the main innovation in this paper. 

\begin{definition}
    A homeomorphism $F:\prod_{i=1}^n \widetilde{X}_{i}\rightarrow \prod_{i=1}^m \widetilde{Y}_{i}$ between a product of metric spaces is said to \emph{preserve the product structure} if up to reindexing the factors there are  homeomorphisms $F_i$ such that for each $i$, the following diagram commutes: \[
    \begin{tikzcd}
\prod_{i=1}^n \widetilde X_{i} \ar[r, "F"] \ar[d,"\pi_i"] & \prod_{i=1}^n \widetilde{Y}_{i} \ar[d,"\pi_i"] \\
\widetilde{X}_{i} \ar["F_i",r] & \widetilde{Y}_{i}.
\end{tikzcd}
\]
\end{definition}

\begin{theoremintro}\label{thm:KKLgeneralizedNEW}
 Let $X_i$ and $Y_j$ be geodesic metric spaces, and $u$ be a subadditive nondecreasing sublinear function. Suppose $\phi:X = \prod_{i=1}^n X_i\rightarrow \prod_{i=1}^n Y_j = Y$ is an $O(u)$-bilipschitz equivalence such that for all $u$-admissible cones $X_\omega$ and $Y_\omega$, the cone map $\operatorname{Cone}(\phi):X_\omega\rightarrow Y_\omega$ preserves the product structure. 
Then there exists a bijection $\sigma: \{ 1, \ldots, n \} \to \{ 1, \ldots, n \}$ and, for every $i \in \{1,\ldots, n\}$, an $O(u)$-bilipschitz equivalence $\phi_i: X_i \to Y_{\sigma(i)}$ such that the diagram
\[ \begin{tikzcd}
 \prod_{i=1}^n X_i \ar[r, "\phi"] \ar[d,"\pi_i"] & \prod_{i=1}^n Y_i \ar[d,"\pi_{\sigma(i)}"] \\
X_i \ar["\phi_i",r] & Y_{\sigma (i)}
\end{tikzcd} \]
commutes up to sublinear error $v$, that is subadditive and nondecreasing. Especially, $\phi$ is $O(v)$-close to the product map $\Phi = \prod_{i=1}^n \phi_i$ of Definition~\ref{def:commuting-sublinear-error}.
\end{theoremintro}

Most of Section~\ref{sec: product theorem} is devoted to the proof of Theorem~\ref{thm:KKLgeneralizedNEW}. 
The other ingredient needed for Theorem~\ref{thm:KKLgeneralized} is the \emph{topological splitting} result from~\cite{KKL}, stating that the spaces considered in Theorem~\ref{thm:KKLgeneralized} satisfy the 
hypothesis of Theorem~\ref{thm:KKLgeneralizedNEW}.

\begin{theorem}[\cite{KKL}, Theorem 5.1]\label{thm: topological splitting}
    Suppose $X_i$ and $Y_j$ are spaces of coarse type I or II. Denote $\overline{X}:=\prod_{i=1}^n X_i, \overline{Y}:=\prod_{i=1}^n Y_i$, and let $\overline{X}_\omega,\overline{Y}_\omega$ be respective asymptotic cones. Suppose $F:=\mathbf{R}^p\times\overline{X}_\omega\rightarrow \mathbf{R}^q\times\overline{Y}_\omega$ is a homeomorphism. Then $p=q$, $n=m$, and the homeomorphism $\overline{F}:\overline{X}_\omega\rightarrow \overline{Y}_\omega$ given by $\pi_{\overline{Y}_\omega}\circ F_{\restriction \overline{X}_\omega}$ preserves the product structure. 
\end{theorem}

Theorem~\ref{thm:KKLgeneralized} is proved in Section~\ref{sec: proving theorem KKL}. In addition to Theorem~\ref{thm:KKLgeneralizedNEW} and Theorem~\ref{thm: topological splitting}, it uses  Proposition~\ref{prop: product with Euclidean factors}, which allows one to deal with Euclidean factors.
\begin{remark}
    In the statements above, we do not explicitly specify how the distance on the product is built from that of the factors. In the proof we work with the $\ell^2$ (or Pythagorean product) metric, which is more natural when the factors are Riemannian. However, the theorem applies if one considers any distances on $X$ and $Y$ quasiisometric to the $\ell^2$ distances, e.g.\ the $\ell^1$ distance.
\end{remark}

The following remarks concern the two differences between our product theorem and the one for quasiisometries in~\cite{KKL}.

\begin{remark}\label{remark: O(1) vs O(u) error term}
    When $\phi$ is a quasiisometry, i.e.\ an $O(1)$-equivalence, the diagrams in the statements of Theorems~\ref{thm:KKLgeneralized} and~\ref{thm:KKLgeneralizedNEW} commute up to a bounded error $D\in O(1)$ (\cite{KKL}, Theorem~A). For $O(u)$-equivalences we can only bound this error term by an arbitrary sublinear function $v$ that may not be in $O(u)$. We use Lemmas~\ref{lem:dominating-sublinear-by-concave} and~\ref{lem:dominating-sublinear-by-subadditive} to upgrade $v$ to a nondecreasing subadditive function, which are natural assumptions in the context of sublinear bilipschitz equivalence~\cite{pallier2018large,GrayevskyRigidity}. We are not aware of any application that may be affected by the fact that $v\notin O(u)$, not even when $u=1$ and $\phi$ is a quasiisometry. We stress that nonetheless, the factor maps $\phi_i:X_i\rightarrow Y_i$ are $O(u)$-equivalences. 
\end{remark}

\begin{remark}\label{rmk: nontranlatability}
One notable and delicate difference between our statements those of \cite{KKL} concern the \emph{nontranslatability} condition: a pair of metric spaces $(Z,W)$ is \emph{nontranslatable} if any two homeomorphisms $f,g:Z\rightarrow W$ that are at finite distance must coincide (\cite[Definitions~2.2, 2.3]{KKL}). In~\cite{KKL}, this notion  is used in two separate ways. First, it is one of the key definitions in the context of their topological splitting (Theorem~\ref{thm: topological splitting} above). The second use of nontranslatability in~\cite{KKL} is as a means to upgrade the error term in the commuting diagrams of Theorems~\ref{thm:KKLgeneralized} and~\ref{thm:KKLgeneralizedNEW} to an $O(1)$ error. 

In our work we do not require this notion of nontranslatability outside of Theorem~\ref{thm: topological splitting}, which may be taken as a black box in this paper. While the assumptions of Theorem~\ref{thm:KKLgeneralizedNEW} must be strengthened in order to ensure that the error term is in $O(u)$ (see Example~\ref{ex: QI of R2}), we could not determine whether nontranslatability of the pairs of factors would be a sufficient assumption.

\end{remark}

\begin{example}\label{ex: QI of R2}
Let $\mathbf{R}_1,\mathbf{R}_2$ be two copies of the line $\mathbf{R}$ (a non non-translatable space). Consider the map $\phi:\mathbf{R}_1\times \mathbf{R}_2\rightarrow \mathbf{R}_1\times \mathbf{R}_2$ given by $(x,y)\mapsto (x+y^{\frac{1}{3}}, y+x^{\frac{1}{3}})$. One may check it is a quasiisometry, and that in any asymptotic cone, the cone map preserves the product structure. The map $\phi_1:\mathbf{R}_1\rightarrow \mathbf{R}_1$ defined as $\phi_1(x)=\pi_1\circ\phi(x,0)$ is a quasiisometry. However, for all $y \in \mathbf R$
$$d(\phi_1\circ\pi_1(x,y),\pi_1\circ\phi(x,y))= \vert y\vert ^{\frac{1}{3}}$$
which is not a bounded function of $\vert y \vert$.

\end{example}

Before stating the second application of Theorem~\ref{thm:KKLgeneralized}, we need to recall Cornulier's $\rho_1$ reduction (see \S\ref{sec:cornulier-reduction-partial-converse} for a  more comprehensive account). Given a completely solvable group $S$ we denote by $\operatorname{R}_{\exp} S$ the smallest normal subgroup such that $S/\operatorname{R_{\exp} }S$ is nilpotent; this is also the intersection of the descending central series of $S$. The homomorphism $\alpha\colon S \to \operatorname{Aut}(\operatorname{R}_{\exp} S) = \operatorname{Aut}(\operatorname{Lie}(\operatorname{R}_{\exp} S))$ determined by the Adjoint action of $S$ is algebraic over $\mathbf R$, and can be decomposed into $\alpha = \alpha_\sigma  \alpha_\nu$, where $\alpha_\sigma$ is valued in a diagonal $\mathbf R$-torus of $\operatorname{Aut}(\operatorname{Lie}(\operatorname{R}_{\exp} S))$, and $\alpha_\nu$ is valued in the unipotent radical of $\operatorname{Aut}(\operatorname{Lie}(\operatorname{R}_{\exp} S))$.
Since $\operatorname{R}_{\exp} S$ is nilpotent (it is contained in $[S,S]$), it sits in $\ker \alpha_\sigma$, so that $\alpha_\sigma$ defines a homomorphism $S /\operatorname{R}_{\exp} S \to \operatorname{Aut}(\operatorname{R}_{\exp} S)$, that we still denote $\alpha_\sigma$.

\begin{definition}
    Let $S$, $\alpha_\sigma$ and $\alpha_\nu$ be as above. One defines $\rho_1(S)$ as the group $\operatorname{R}_{\exp} S \rtimes_{\alpha_{\sigma}} S/\operatorname{R}_{\exp} S$. We say that $S$ is in the class $(\mathcal C_1)$ if $S\simeq \rho_1(S)$, that is, if $\operatorname{R}_{\exp} S$ is split and $\alpha_\nu = 1$. (In \cite{CornulierCones11}, $\rho_1(S)$ is denoted $S_1$.)
\end{definition}

\begin{theorem}[{Cornulier, \cite{CornulierCones11}}]
\label{th:Cornulier-thm-intro}
Let $S$ be a completely solvable group.
Then $S$ and $\rho_1(S)$ are $O(\log)$-bilipschitz equivalent.
\end{theorem}

Combining this theorem with our product theorem, and using some further previous results on sublinear bilipschitz equivalences, we obtain the following.

\begin{theoremintro}\label{cor:KKLgeneralized}
    Let $S$ and $S'$ be two completely solvable groups.
    Assume that there exists a sublinear, subadditive function $u$ such that $S$ and $S'$ are $O(u)$-sublinear bilipschitz equivalent (which holds, in particular, if $S$ and $S'$ are quasiisometric with $u=1$ in that case).
    Further assume that
    \[ \rho_1(S) \simeq \mathbf R^n \times P \times H_1 \times \cdots \times H_m \quad \text{and} \quad\rho_1(S') \simeq  \mathbf R^{n'} \times P' \times H'_1 \times \cdots \times H'_{m'} \]
    for some $n,m,n',m'\geqslant 0$, where
    \begin{enumerate}
        \item \label{item:assum1} $P=AN$ and $P'=A'N'$ are maximal completely solvable subgroups in semisimple groups $G=KAN$ and $G'=K'A'N'$ respectively.
        \item \label{item:assum2} For $i=1,\ldots, m$, $H_i$ has a left-invariant Riemannian metric that is negatively curved, and an abelian derived subgroup; same assumption for $H'_j$, $j = 1,\ldots, m'$.
    \end{enumerate}
    Then, $\rho_1(S)$ and $\rho_1(S')$  are isomorphic.
\end{theoremintro}

The case $S=P$ is equivalent to the former Corollary~\ref{cor:KKL-symmetric} while the case where $S$ is equal to a single factor $H_1$ is the main theorem of \cite{pallier2019conf}, used in the proof.

A more sophisticated (and slightly more general) version of Theorem \ref{cor:KKLgeneralized} will be given in Theorem \ref{th:technical-KKL}; in particular, the geometric assumption on curvature in \eqref{item:assum2} can be reformulated in terms of the structure of $H$ with the help of Heintze's theorem \cite{Heintze}. 

Theorem~\ref{thm:KKLgeneralized} and Theorem~\ref{cor:KKLgeneralized} allow us to distinguish between several families of completely solvable groups up to quasiisometry. The following is an example of this strategy in dimension $4$.

\begin{example}\label{exm:four-dim-KKL}
Let $\alpha \in (0,1)$.
    The groups $S=G^0_{4,9}$ and $S'_{\alpha}=\mathbf R \times G_{3,5}^{\alpha}$ (the names are from the classification in \cite{Mubarakzyanov}) are the four dimensional, completely solvable Lie groups whose respective Lie algebras $\mathfrak g_{4,9}^0$ and $\mathbf R \times \mathfrak g_{3,5}^{\alpha}$ are spanned by $e_1, \ldots, e_4$, subject to the following nonzero brackets:
    \begin{align*}
        \mathfrak g_{4,9}^0: \quad 
        & [e_4,e_1] = e_1, [e_4,e_2] = e_2, [e_2, e_3] = e_1 \\
        \mathbf R \times \mathfrak g_{3,5}^{\alpha}:\quad  & [e_4,e_1] = e_1, [e_4,e_2] = \alpha e_2.
    \end{align*}
    It follows from Theorem \ref{cor:KKLgeneralized} (with $u \equiv 1$) that for every $\alpha \in (0,1)$, $S$ and $S'_{\alpha}$ are not quasiisometric. Here, $\rho_1(S)$ splits as a direct product, that is $\mathbf R \times P$ where $P$ is the maximal completely solvable subgroup of $\mathrm{SO}(3,1)$. In this example the original Kapovich-Kleiner-Leeb theorem cannot be applied since $S$ is not a direct product of non-trivial Lie groups. See section~\ref{sec: example 1.3} for a direct proof of this example, that we use as a warm-up for the general case.
\end{example}

\begin{remark}
The parameter $\alpha$ in the definition of $G_{3,5}^{\alpha}$ can actually be taken in $(-1,1) \setminus \{ 0 \}$; the two limit cases $\alpha = -1$ and $\alpha =1$ give the groups $\mathbf R \times \mathrm{Sol}_3$ and $\mathbf R \times P$ respectively, where $P$ is as above.
A direct inspection of the asymptotic cone suffice to tell $G_{4,9}^0$ apart from the groups $\mathbf R \times G_{3,5}^{\alpha}$ when $\alpha <0$ and from $\mathbf R\times \mathrm{Sol}_3$.
As the group $\mathbf R \times P$ bears a symmetric metric with a Euclidean factor, one can describe the finitely generated groups quasiisometric to $\mathbf R \times P$ \cite{KleinerLeebCAG}; however, due to the presence of the Euclidean factor we do not know whether similar methods can describe the Lie groups quasiisometric to $\mathbf R \times P$; we raise this question in Appendix \ref{app:completely-solvable-symmetric}.
\end{remark}

For completely solvable groups, primary quasiisometry invariants are the cone dimension (that is, the covering dimension of the asymptotic cone) and the Dehn function. In a subsequent paper~\cite{grayevsky2025dehnfunctionscomputationslower}, we use the wide range of tools developed in~\cite{CoTesDehn} by Cornulier and Tessera in order to determine the Dehn function of all completely solvable groups of exponential growth and dimension up to 5. In Corollary~\ref{cor: contribution of product theorem} below we will summarize the contribution of our strategy to the quasiisometric classification of these groups. More precisely, Corollary~\ref{cor: contribution of product theorem} lists the $5$-dimensional groups which share the same cone dimension and Dehn function, but for which Theorem \ref{cor:KKLgeneralized} implies that they are nonetheless quasiisometrically distinct. We recover a fact already obtained recently (in a more specific form that we explain below) by Bourdon and Rémy {\cite{bourdon2023rhamlpcohomologyhigherrank}}.

\begin{corollaryintro}\label{corintro: contribution of product theorem} 
    There exists uncountably many quasiisometry classes of indecomposable, non-unimodular, completely solvable Lie groups with quadratic Dehn function.
\end{corollaryintro}

The family of pairwise non-isomorphic groups we consider in order to deduce this corollary are named $G_{5,19}^{1,\beta}$ with $\beta\in (0,+\infty)$ in \cite{Mubarakzyanov} while Bourdon and Rémy's family is $G_{5,33}^{\alpha, 1-\alpha}$ with $\alpha \in [1/2, +\infty)$ (it is indecomposable for $\alpha\neq 1$).
All these groups have dimension 5 and cone dimension $2$, but the groups considered by Bourdon and Rémy have left-invariant proper CAT(0) metrics (a condition known to imply a quadratic isoperimetric inequality, see e.g. \cite{WengerIsopEuclid}), while the groups in the $G_{5,19}^{1,\beta}$ family do not have such metrics, as can be checked using Azencott and Wilson's criterion \cite{azencott1976homogeneous}.
Note that Peng obtained uncountably many quasiisometry classes among the completely solvable groups \cite{PengCoarseI,PengCoarseII}.
The groups in the $G_{5,19}^{1,\beta}$ family fall outside of the CAT(0) groups and those considered by Peng, which are unimodular.

Bourdon and Rémy obtain their result through the determination of a critical exponent for $L^p$-cohomology in degree $2$; by using \cite{pallier2019conf}, our proof also ultimately involves (although not in the present paper) the use of a critical exponent. The latter appears at first in \cite{pallier2019conf} as a variant of Pansu's conformal dimension; however it coincides with the critical exponent for $L^p$-cohomology in degree one of a family Heintze groups occuring as direct factors of $\rho_1(G_{5,19}^{1,\beta})$ (See  \cite{BourdonKleinerCLPi}, \cite{CarrascoOrliczHeintze} and \cite[\S 3.2.1]{pallier2019conf}). It would thus be interesting to investigate the $L^p$-cohomology of $G_{5,19}^{1,\beta}$ in degree $2$.

In order to formulate our final question, let us call a group of type NPC a completely solvable group with a CAT(0) left-invariant proper Riemannian metric, as characterized in \cite{azencott1976homogeneous}. 
Examples of groups of type NPC are provided by the $AN$ subgroups of (high rank) semisimple Lie groups $G=KAN$, but this family is considerably more vast; in dimension $5$ already, it includes the uncountable family studied by Bourdon and Rémy and mentionned above.

\begin{ques}\label{ques:aw-factors}
    In Theorem~\ref{thm:KKLgeneralized}, can one allow the factors to be taken among non-abelian, indecomposable groups of type NPC?
\end{ques}

The question can equally be raised in the weaker form, for the original version of the Kapovich-Kleiner-Leeb product theorem instead of Theorem~\ref{thm:KKLgeneralized}.

\begin{figure}
    \centering
    \begin{tikzpicture}[line cap=round,line join=round,>=angle 45,x=1.0cm,y=1.0cm]
\clip(-5.5,-4.1) rectangle (9.5,1.9);
\draw (-0.9,1.29) node[stuff_fill,anchor=north east] {\begin{minipage}{4cm}{Topological splitting, Theorem~\ref{thm: topological splitting} \cite{KKL}} \end{minipage}};
\draw (2.06,0.79) node[stuff_fill,anchor=north] {Theorem~\ref{thm:KKLgeneralizedNEW} };
\draw (5.06,1.29) node[stuff_fill,anchor=north west] {\begin{minipage}{3.5cm}{$\rho_1$-- reduction, Theorem~\ref{th:Cornulier-thm-intro} \cite{CornulierCones11} } \end{minipage}};
\draw (0.04,-1.07) node[stuff_fill,anchor=center] {Theorem~\ref{thm:KKLgeneralized}    };
\draw (4.24,-1.07) node[stuff_fill,anchor=center] {Theorem~\ref{th:technical-KKL} };
\draw (4.04,-2.27) node[stuff_fill,anchor=center] {Theorem~\ref{cor:KKLgeneralized} };
\draw (0.04,-2.58) node[stuff_fill,anchor=north] {\begin{minipage}{4.5cm}{SBE classification of symmetric spaces, Corollary~\ref{cor:KKL-symmetric}} \end{minipage}};
\draw (4.04,-3.2) node[stuff_fill,anchor=north] {Corollary~\ref{corintro: contribution of product theorem} };
\draw [->] (-3,0) -- (-1.2,-1);
\draw [->] (2,0) -- (0.5,-0.5);
\draw [->] (1.2,-1.1) -- (2.9,-1.1);
\draw [->] (6.5,0) -- (5.6,-1);
\draw [->] (0.1,-1.5) -- (0.1,-2.5);
\draw [->] (4.5,-1.5) -- (4.5,-1.9);
\draw [->] (4.5,-2.7) -- (4.5,-3.1);
\draw [dash pattern = on 2pt off 2pt] (2.6,-0.5) -- (6,-0.5) -- (6,-4) -- (2.6,-4) -- (2.6,-0.5);
\draw (5.5,-3.1) node[stuff_fill_blank,anchor=west] {\begin{minipage}{2.5cm}{Applications to solvable Lie groups} \end{minipage}};
\draw (-5,-2) node[stuff_fill_blank,anchor=west] {\begin{minipage}{4.5cm}{ with \cite{pallier2018large} (rank one) \\ and \cite{GrayevskyRigidity} (no rank one)} \end{minipage}};
\end{tikzpicture}
    \caption{Logical dependencies between the main statements in this paper.}
    \label{fig:placeholder}
\end{figure}

\subsection{Organization of the paper}
The proofs of Theorems~\ref{thm:KKLgeneralized} and~\ref{thm:KKLgeneralizedNEW} are in Section~\ref{sec: product theorem}. Section~\ref{sec: completely solvable groups and corollaries} introduces some of the theory of completely solvable Lie groups and the significance of sublinear bilipschitz equivalence to this theory. Corollary~\ref{cor:KKL-symmetric} and Theorem~\ref{cor:KKLgeneralized} are proved in Section~\ref{subsec:proof-corollaries}. We actually deduce both of these statements from Theorem~\ref{th:technical-KKL}, although Corollary~\ref{cor:KKL-symmetric} could be derived directly. Theorem~\ref{th:technical-KKL} is stated in terms of diagonal Heintze groups, which are the Gromov-hyperbolic groups of class $(\mathcal C_1)$ and thus are of coarse type I. Finally, in Section~\ref{sec: contribution of product theorem} we formulate and prove Corollary~\ref{cor: contribution of product theorem} about the families of $5$-dimensional Lie groups that can be distinguished up to quasiisometry using our product theorem; Corollary \ref{corintro: contribution of product theorem} appears as a byproduct of this study, in Corollary~\ref{cor: g519beta}.

\subsection{Acknowledgments}
We thank Yves Cornulier for a useful discussion. The second-named author thanks Marc Bourdon and Bertrand Rémy for sharing their results, and Antoine Velut for useful discussions.

\section{The product Theorem \ref{thm:KKLgeneralized}}\label{sec: product theorem}

In this section we prove Theorem~\ref{thm:KKLgeneralizedNEW} and conclude Theorem~\ref{thm:KKLgeneralized}. We start with some preliminaries on sublinear bilipschitz equivalence and sublinear functions. Then we present asymptotic cones, which are the key objects in the proof. We then turn to prove Theorem~\ref{thm:KKLgeneralizedNEW} and formulate (without proof) a variant which allows us to deal with Euclidean factors in Proposition~\ref{prop: product with Euclidean factors}. In Section~\ref{sec: proving theorem KKL} we deduce Theorem~\ref{thm:KKLgeneralized}. The strategy appears already in~\cite{KKL}, and most of our work is in adapting it to the sublinear setting. 

\subsection{Preliminaries}
Let $u \colon \mathbf R_{\geqslant 0} \to \mathbf R_{\geqslant 1}$ be a sublinear function, that is, 
\[ \lim_{r \to + \infty} \frac{u(r)}{r} =0. \]

We say that $u$ is admissible if it is nondecreasing and $\limsup u(2r)/u(r)<+\infty$.
Often, we will also assume that $u$ is subadditive, that is, $u(r_1+r_2) \leqslant u(r_1) + u(r_2)$ for all $r_1,r_2 \geqslant 0$. This implies admissibility.
This assumption is not extremely restrictive, as we prove in Lemma~\ref{lem:dominating-sublinear-by-subadditive} below: a sublinear nondecreasing function is always dominated by a subadditive (in particular, admissible) one.
We first need an intermediate technical statement.

\begin{lemma}
\label{lem:dominating-sublinear-by-concave}
    Let $u_0 \colon \mathbf R_{\geqslant 0} \to \mathbf R_{\geqslant 0}$ be such that $\lim_{r \to +\infty} \frac{u_0(r)}{r} = 0$, $u_0(0)=0$, and $\limsup_{r \to 0} \frac{u_0(r)}{r} < +\infty$. Then there exists $v_0 \colon \mathbf R_{\geqslant 0} \to \mathbf R_{\geqslant 0}$ such that $u_0 \leqslant v_0$ and $v_0$ is a concave function such that $\lim_{r \to +\infty} \frac{v(r)}{r} = 0$. 
\end{lemma}

\begin{proof}
    Let $\mathcal W$ be the set of concave functions $w$ from $\mathbf R_{\geqslant 0}$ to itself and such that $u_0 \leqslant w$. Given the assumptions that we made on $u_0$, $\mathcal W$ contains a linear function, especially it is not empty.
    For all $r\in \mathbf R_{\geqslant 0}$ we set $v_0(r) = \inf \{ w(r) \colon w \in  \mathcal W\}$ and proceed to prove that the function $v_0$ is sublinear. If it was not the case, there would be $\varepsilon >0$ and a sequence $\{ r_i \}_{i \in \mathbf N}$ going to infinity, such that $v_0(r_i) \geqslant \varepsilon r_i$. Note however that $v_0$ is concave, as it is an infimum of concave functions. So from the previous inequality and the convexity inequality, we get that $v_0 (r) \geqslant r$ for every nonnegative $r$. 
    But $u_0$ is sublinear, so there exists $R > 0$ such that $u_0(r) \leqslant \varepsilon r /2$ for every $r \geqslant R$. 
    Considering the function $v_1$ such that $v_1(r) = \epsilon r$ if $0 \leqslant r \leqslant R$ and $v_1(r) = \varepsilon R/2 +  \varepsilon r/2 $ if $r \geqslant R$, we get a function $v_1$ in $\mathcal W$ eventually smaller than $v_0$, a contradiction. This terminates the proof that $v_0$ is sublinear; the other required properties of $v_0$ have been established in the course of the proof.
\end{proof}

\begin{lemma}[Dominating sublinear increasing functions by subadditive ones]
\label{lem:dominating-sublinear-by-subadditive}
    Let $u \colon \mathbf R_{\geqslant 0} \to \mathbf R_{\geqslant 1}$ be a nondecreasing sublinear function. There exists a subadditive, nondecreasing, sublinear function $v$ such that $u \leqslant v$.
\end{lemma}

\begin{proof}
    Given the assumptions made on $u$, there exists $M \geqslant 0$ such that $u-M \leqslant u_0$, where $u_0$ is a sublinear function such that $u_0(0) = 0$ and $\limsup_{r \to 0} \frac{u_0(r)}{r} = 0$.
    Applying Lemma~\ref{lem:dominating-sublinear-by-concave} to $u_0$ we get a concave function $v_0$ dominating $u_0$. Now set $v_1 = v_0+M$. We claim that $v_1$ is subadditive. Indeed, $v_0$ is concave and $v_0(0)=0$, so that for every pair $(r,s)$ of nonnegative real numbers, not both equal to zero, we have
\[ v_0(r) \geqslant \frac{r}{r+s}v_0(r+s) \; \text{and}\; v_0(s) \geqslant \frac{s}{r+s}v_0(r+s) \]    
which after summation gives that $v_0$ is subadditive; now, for any $r,s \geqslant 0$,
    \begin{align*}
        v_1(r+s) & = v_0(r+s) + M \\
        & \leqslant v_0(r) + v_0 (s) + M \\
        & \leqslant v_0 (r) + M + v_0(s) + M \\
        & = v_1(r) + v_1(s).
    \end{align*}
    In addition, $v_0$ is sublinear, so $v_1$ is as well. Finally, assume towards contradiction that there exists $(r,s)$ such that $0 \leqslant r < s$ and $v_1(r) > v_1(s)$. $v$ being concave, it will then eventually take negative values, as $v_1(t) \leqslant v_1(r) -\frac{t-r}{s-r}(v_1(r)-v_1(s))$ which becomes negative for large enough $t$. This is a contradiction, hence $v_1$ is nondecreasing. One may take $v = v_1$ and the conclusion is reached.
\end{proof}

\begin{remark}
    We do not have $v \in O(u)$ in the conclusion of lemma~\ref{lem:dominating-sublinear-by-subadditive}, as the following example shows: $u(t) = \exp(\lfloor \log (\log (t))\rfloor^2)$ for large enough $t$. Then for $t_n = e^{e^n} -1/2$ we have $\lim_n u(t_n+1)/u(t_n) = +\infty$, and $\lim_n v(t_n+1)/v(t_n) = +\infty$ for any $v \in O(u)$ dominating $u$, which cannot be for a subadditive $v$.
\end{remark}

\begin{remark}
    The class of admissible functions is stable under the operations $\vee$ and $\wedge$. This will be of use in Section~\ref{sec: completely solvable groups and corollaries} where we will use that if $u$ is admissible, then $u \vee \log$ is admissible.
\end{remark}

Let $X$ and $Y$ be  metric spaces. After fixing $x_0 \in X$ and $y_0 \in Y$, we denote by $\vert \cdot \vert$ the distance to the respective basepoints in $X$ and $Y$. We denote $\vert x_1\vert\vee \vert x_2\vert:=\max\{\vert x_1\vert,\vert x_2\vert\}$.

\begin{definition}\label{def: SBE}
Let $X$, $Y$, $x_0$, $y_0$ and $u$ be as above. Let $L \geqslant 1$.
    We say that $f\colon X \to Y$ is 
    \begin{itemize}
        \item $(L,u,x_0)$-Lipschitz if for every $x,x' \in X$,
    $d\big(f(x),f(x')\big) \leqslant Ld(x,x') + u (\vert x \vert \vee \vert x' \vert)   $
    \item $(L,u,x_0)$-expansive 
    if $L^{-1}d(x,x') - u (\vert x \vert \vee \vert x' \vert) \leqslant d\big(f(x), f(x')\big)$
    \item $(u,y_0)$-surjective if for every $y$ in $Y$, there is 
    $x \in X$ such that $d\big(y,f(x)\big) \leqslant u(\vert y \vert)$.
    \end{itemize}

\end{definition}

We say that $f$ is an $(L,u)$-bilipschitz embedding if it is $(L,cu,x_0)$-Lipschitz and $(L,cu,x_0)$-expansive for some $c\geqslant 0$.
If $f$ is additionally $(u,y_0)$-surjective for some $y_0$, then for all $y'_0 \in Y$ there is $c'>0$ such that it is $(c'u, y'_0)$-surjective; in this case, we say that $f$ realizes an $\big(L,O(u)\big)$-bilipschitz equivalence, or for short, an $O(u)$-bilipschitz equivalence between $X$ and $Y$. In line with the terminology for quasiisometries, $L$ is sometimes called the bilipschitz constant of $f$.
When no reference is made to $L$ and $u$ we will call an $(L,u)$-bilipschitz embedding a sublinear bilipschitz embedding. We may omit the base points or any of the constants when they are clear from the context.

\subsection{Going through cones}

\begin{definition}\label{def: cone}
    Let $X$ be a metric space. Let  $(\sigma_k)_{k \in \mathbf N}$ be a sequence of positive real numbers. Let $(x_k)_k \in X^{\mathbf N}$.
    We call precone of $X$ with data $\big((x_k)_k, (\sigma_k)_k\big)$ the set of sequences
    \[ \operatorname{Precone}\big(X,(x_k)_k, (\sigma_k)_k\big) = \{ (x'_k)_k \in X^{\mathbf N} : \exists M \in [0,+\infty), \forall k \in \mathbf N, d(x_k,x_k') \leqslant M\sigma_k \}. \]
    Given a nonprincipal ultrafilter $\omega$ over $\mathbf N$ we denote by $\operatorname{Cone}\big(X, (x_k)_k, (\sigma_k)_k,\omega\big)$ the quotient of the set $\operatorname{Precone}\big(X,(x_k)_k, (\sigma_k)_k\big) $ by the relation 
    \[ (x'_k)_k \sim (x''_k)_k \iff \lim_\omega \frac{d(x'_k,x''_k)}{\sigma_k} = 0 \]
    equipped with the distance 
    \[ d_\omega(\mathbf [x'], \mathbf [x'']) =  \lim_\omega \frac{d(x'_k,x''_k)}{\sigma_k}. \]
\end{definition}
If $\sigma_k$ goes to $+\infty$ we call this an asymptotic cone.

When it is not relevant to mention all the remaining data, we will denote an asymptotic cone of $X$ by $X_\omega$. Often the defining data is relevant, in which case we ease notation and write $\operatorname{Cone}(X,x_k,\sigma_k)$.

\begin{definition}\label{def: u admissible}
    Given a metric space $X$, a triple $(X,x_k,\sigma_k)$ is called $u$-admissible if for some (any) $v\in O(u)$ 
    \begin{equation*}
    \lim_{k \to + \infty} \frac{v(|x_k|)}{\sigma_k} = 0.
\end{equation*}
An asymptotic cone $X_\omega:=\operatorname{Cone}(X,x_k,\sigma_k)$ is called $u$-admissible if the triple $(X,x_k,\sigma_k)$ is $u$-admissible. 
\end{definition}

Note that we did not specify a basepoint when writing $v(|x_k|)$ in the above definition; such a choice turns out to have no influence on the notion of $u$-admissibility.

Lemma~\ref{lem:go-through-cones} below is a basic fact, usually applied to cones with fixed base point. We will need asymptotic cones with moving base points, which are slightly less common.
In the statement below, given a map $f:X \to Y$ we still denote $f$ the map between the power sets $f: X^\mathbf N \to Y^{\mathbf N}$.

\begin{lemma}
\label{lem:go-through-cones}
Let $L \geqslant 1$.
Let $u$ be an admissible function, let $X$ and $Y$ be metric spaces and let $f:X\rightarrow Y$ an $(L,u)$-bilipschitz equivalence.
Assume $\big(X,(x_k)_k,(\sigma_k)_k\big)$ is $u$-admissible. Then
\begin{equation*}
    f \Big(\operatorname{Precone}\big(X,(x_k)_k,(\sigma_k)_k\big)\Big) \subseteq \operatorname{Precone}\Big(Y,\big(f(x_k)\big)_k, (\sigma_k)_k\Big),
\end{equation*}
and for any nonprincipal ultrafilter $\omega$, a quotient map 
\[\operatorname{Cone}\big(f,(x_k)_k, (\sigma_k)_k, \omega\big) \colon \operatorname{Cone}\big(X,(x_k)_k, (\sigma_k)_k,\omega) \to \operatorname{Cone}\Big(Y,\big(f(x_k)\big)_k, (\sigma_k)_k, \omega\Big) \]
is well-defined and an $L$-bilipschitz embedding.
Moreover, if $f$ is additionally $O(u)$-surjective, then $\operatorname{Cone}\big(f,(x_k)_k, (\sigma_k)_k\big)$ is a bilipschitz homeomorphism.
\end{lemma}

\begin{proof}
    This follows from \cite[Proposition 6]{loglie}.
\end{proof}

\begin{remark}\label{importance of u admissibility to the error term}
    The notion of $u$-admissibility concerns a key difference between quasiisometries and $O(u)$-bilipschitz equivalence. When $u\in O(1)$, then every asymptotic cone is $u$-admissible, and in particular the cone maps of quasiisometries are bilipschitz homeomorphisms in every asymptotic cone. However when $u\notin O(1)$, not all cone maps are bilipschitz homeomorphisms. In our work, this adds a significant difficulty in the context of controlling the error term in Theorem~\ref{thm:KKLgeneralizedNEW} (see Remarks~\ref{remark: O(1) vs O(u) error term} and~\ref{rmk: nontranlatability}).
\end{remark}

\begin{remark}\label{rem: cone notations}
    Unless stated otherwise, whenever we define $X_\omega=\operatorname{Cone}(X,x_k,\sigma_k)$, we implicitly define $Y_\omega:=\operatorname{Cone}\big(Y,f(x_k),\sigma_k\big)$. Lemma~\ref{lem:go-through-cones} implies that if $X_\omega$ is $u$-admissible, so is $Y_\omega$. 
\end{remark}

\subsection{Proof of Theorem~\ref{thm:KKLgeneralizedNEW}} \label{sec: generalizing the proof of KKL}

Throughout this section, the standing assumptions on the map $\phi$ and spaces $X,Y$ are those of Theorem~\ref{thm:KKLgeneralizedNEW}. 
When it is clear from the context, we occasionally do not specify the space in which the distance is computed. Also, we sometimes write $\pi_j$ instead of $\pi_{X_j}$ or $\pi_{Y_j}$. 

The following are some sublinear adaptations to the notions that appear in~\cite{KKL}. 

\begin{definition}
Let $d>0$ be a positive constant and $d_0:\mathbf{R}_{\geq 0}\rightarrow \mathbf{R}_{\geq 1}$ a function. A pair of points $x,x'\in X$ is called \emph{$d$-separated} if $d_X(x,x')\geqslant d$, and \emph{$d_0$-separated} if $d(x,x')\geqslant d_0(|x|\vee|x'|)$.
\end{definition}

\begin{definition}
    A pair of points $x,x'\in X=\prod_i X_i$ is called \emph{$i$-horizontal} if $x$ and $x'$ agree on all $m\ne i$ coordinates.
\end{definition}
 
\begin{definition}\label{def: compressability}
    Let $x,x'\in X$ be a pair of $i$-horizontal points, and fix $L_1>1, \epsilon<L_1^{-1}$. The following definitions depend on $L_1$ and $\epsilon$, but we suppress them from the terminology:
    
    \begin{enumerate}
        \item $(x,x')$ is \emph{$j$-compressed} if 

    $$\frac{d_{Y_j}\Big(\pi_{Y_j}\big(\phi(x)\big),\pi_{Y_j}\big(\phi(x')\big)\Big)}{d_X(x,x')}\leq \epsilon $$

    \item $(x,x')$ is  \emph{$j$-uncompressed} if
    
    $$L_1^{-1}\leq \frac{d_{Y_j}\Big(\pi_{Y_j}\big(\phi(x)\big),\pi_{Y_j}\big(\phi(x')\big)\Big)}{d_X(x,x')}\leq L_1$$

    \item $(x,x')$ is  \emph{$j$-semi-compressed} if
        $$\epsilon \leq \frac{d_{Y_j}\Big(\pi_{Y_j}\big(\phi(x)\big),\pi_{Y_j}\big(\phi(x')\big)\Big)}{d_X(x,x')}\leq {L_1}^{-1}.$$

    \end{enumerate}
 Two pairs are \emph{$j$-compatible} if they are simultaneously $j$-semi/un/compressed, and \emph{$j$-incompatible} otherwise.
\end{definition}

Recall that $L$ is the bilipschitz constant of the map $\phi:X\rightarrow Y$. We fix once and for all $L_1>L$ and $\epsilon<L_1^{-1}$, and use the terminology of Definition~\ref{def: compressability} with respect to these constants. Our proof of Theorem~\ref{thm:KKLgeneralizedNEW} works for any choice of $L_1> L$ and $\epsilon<L_1^{-1}$ (compare below with Remark~\ref{rmk: omitted proof of Euclidean factors proposition}).

For an $i$-horizontal pair $x,x'\in X$, we note two trivial facts:
\begin{enumerate}
    \item  $d_{X_i}\big(\pi_i(x),\pi_i(x')\big)=d_X\big(x,x')$.
    \item Any geodesic $[x,x']$ in $X$ joining $x$ and $x'$ is contained in a fiber of $X_i$. In particular, for each $x''\in [x,x']$, $(x,x'')$ and $(x',x'')$ are $i$-horizontal pairs.
\end{enumerate}

\begin{lemma} \label{epsilon compressed transitivity}
    Being $j$-compressed is transitive along a geodesic, i.e.\ if $x''$ lies on a geodesic joining a pair $(x,x')$ of $i$-horizontal points, and both $(x,x'')$ and $(x'',x')$ are $j$-compressed, then also $(x,x')$ is $j$-compressed.
\end{lemma}

 \begin{proof}
       
\begin{align*}
    d_{Y_j}\Big(\pi_j\big(\phi(x)\big),\pi_j\big(\phi(x')\big)\Big) & \leq d_{Y_j}\Big(\pi_j\big(\phi(x)\big),\pi_j\big(\phi(x'')\big)\Big)+d_{Y_j}\Big(\pi_j\big(\phi(x'')\big),\pi_j\big(\phi(x')\big)\Big) \\ & \leq \varepsilon d(x,x'')+\varepsilon d(x'',x')=\varepsilon d(x,x')
\end{align*}
    
    The last equality follows from $x''\in[x,x']$.
 \end{proof}

The following lemma is straightforward and will be used repeatedly. 

\begin{lemma}\label{lem: cone contradiction}
    Let $X=\prod_{i=1}^n X_i, Y=\prod_{i=1}^n Y_i$, $\phi:X\rightarrow Y$ an $(L,u)$-bilipschitz equivalence. Let $X_\omega=\operatorname{Cone}(X,x_k,\sigma_k), Y_\omega=\operatorname{Cone}(Y,\phi(x_k),\sigma_k)$ be $u$-admissible asymptotic cones, and assume $\phi_\omega =\operatorname{Cone}(\phi):X_\omega\rightarrow Y_\omega$ preserves the product structure, with factor maps $\phi_{\omega_i}:X_{\omega_i}\rightarrow Y_{\omega_i}$. Let $z_\omega:=[(z_k)_k], z_\omega':=[(z_k')_k]\in X_\omega$. It holds that

    \begin{enumerate}
        \item \label{item:first-item-cone-contradiction} If $\lim_\omega \frac{1}{\sigma_k}d_{X_j}\big(\pi_j(z_k), \pi_j(z_k')\big)=0$ then $\lim_\omega \frac{1}{\sigma_k}d_{Y_j}\Big(\pi_j\big(\phi(z_k)\big),\pi_j\big(\phi(z_k')\big)\Big)=0$.

        \item \label{item:second-item-cone-contradiction} If for all $l\ne i$ 
        $\lim_\omega \frac{1}{\sigma_k}d_{X_l}\big(\pi_l(z_k), \pi_l(z_k')\big)=0$, then 
        $$\lim_\omega \frac{d_{Y_i}\Big(\pi_i\big(\phi(z_k)\big),\pi_i\big(\phi(z_k')\big)\Big)}{d_X(z_k,z_k')}\in \left[ \frac{1}{L},L \right] .$$

    \end{enumerate}
   In particular, if the pairs $(z_k,z_k')$ which define $z_\omega$ and $z_\omega'$  are $\omega$-almost surely $i$-horizontal, then $\omega$-almost surely $(z_k,z_k')$ is $j$-compressed for all $j\ne i$ and $i$-uncompressed.
\end{lemma}

\begin{proof}
    In the notations of the statement, the hypothesis $\lim_\omega\frac{1}{\sigma_k}d\big(\pi_i(z_k), \pi_i(z_k')\big)=0$ implies $\pi_{\omega_i}(z_\omega)=\pi_{\omega_i}(z_\omega')$, where $\pi_{\omega_i}$ is the projection to the $i$-th factor in the cone. Therefore the assumption that $\phi_\omega$ preserves the  product structure gives $\pi_{\omega_i}\big(\phi_\omega(z_\omega)\big)=\pi_{\omega_i}\big(\phi_\omega(z_\omega')\big)$. Item \eqref{item:first-item-cone-contradiction} follows. Under the assumption of \eqref{item:second-item-cone-contradiction}, we thus have
    $$d_{Y\omega}\big(\phi_\omega(z_\omega),\phi_\omega(z_\omega')\big)=d_{Y_{\omega_i}}\Big(\pi_{\omega_i}\big(\phi_\omega(z_\omega)\big),\pi_{\omega_i}\big(\phi_\omega(z_\omega')\big)\Big).$$
    Since $X_\omega$ is $u$-admissible,  $\phi_\omega$ is $L$-bilipschitz by Lemma~\ref{lem:go-through-cones}, hence    \begin{equation}\label{eq: distance in the cone of Y_1 is bounded when the other components are equal}
    \frac{1}{L}d_\omega(z_\omega,z_\omega')\leq d_{Y_{\omega_1}}\Big(\pi_{i}\big(\phi(z_\omega)\big),\pi_i\big(\phi(z_\omega')\big)\Big)\leq Ld_\omega(z_\omega,z_\omega'). 
    \end{equation}
    By the definition of the distance in the cone, we have
    \begin{align*}
        \lim_\omega \frac{d_{Y_i}\Big(\pi_i\big(\phi(z_k)\big),\pi_i\big(\phi(z_k')\big)\Big)}{d_X(z_k,z_k')}& = \frac{\lim_\omega\frac{1}{\sigma_k}d_{Y_i}\Big(\pi_i\big(\phi(z_k)\big),\pi_i\big(\phi(z_k')\big)\Big)}{\lim_\omega\frac{1}{\sigma_k}d_X(z_k,z_k')}\\
    & =\frac{d_{Y_{\omega_1}}\Big(\pi_{i}\big(\phi(z_\omega)\big),\pi_i\big(\phi(z_\omega')\big)\Big)}{d_\omega(z_\omega,z_\omega')}\in \left[\frac{1}{L},L \right],
    \end{align*}
    which proves the second item. 
    Assume further that the pairs $(z_k,z_k')$ defining $z_\omega,z_\omega'\in X_\omega$ are $i$-horizontal $\omega$-almost surely. 
    In particular, the assumption of \eqref{item:second-item-cone-contradiction} holds (and automatically the hypothesis of \eqref{item:first-item-cone-contradiction} for $j\ne i$). Since $L_1>L$, this immediately implies that $\omega$-almost surely the pairs $(z_k,z_k')$ are $i$-uncompressed. Fix $j\ne i$. By definition of the cone, $z_\omega,z_\omega'\in X_\omega$ implies $\lim_\omega \frac{d_X(z_k,z_k')}{\sigma_k}=C$ for some $C<\infty$, so there are $\epsilon_k\geq 0$ such that $C\sigma_k-\epsilon_k\leq d_X(z_k,z_k')$ and  $\lim_\omega \frac{\epsilon_k}{\sigma_k}=0$. Therefore $\omega$-almost surely 
    $$\sigma_k\frac{C}{2}<\sigma_k(C-\frac{\epsilon_k}{\sigma_k})=\sigma_kC-\epsilon_k\leq d_X(z_k,z_k').$$ 
    We then have that $\omega$-almost surely 
    $$\frac{d_{Y_j}\Big(\pi_j\big(\phi(z_k)\big),\pi_j\big(\phi(z_k')\big)\Big)}{d_X(z_k,z_k')}\leq \frac{d_{Y_j}\Big(\pi_j\big(\phi(z_k)\big),\pi_j\big(\phi(z_k')\big)\Big)}{\sigma_k\frac{C}{2}}.$$
    The ultralimit of the right-hand side is $0$ by the conclusion of the first item of the lemma. 
    Therefore for every $\epsilon$, $\omega$-almost surely the left hand side is smaller than $\epsilon$. We conclude that  $\omega$-almost surely $(z_k,z_k')$ is $j$-compressed for all $j\ne i$.
\end{proof}
The following lemma is the core of the proof of Theorem~\ref{thm:KKLgeneralizedNEW}.  
\begin{lemma}[$O(u)$ version of Lemma~$2.7$ in~\cite{KKL}]\label{lem: simultaneously-compressed-decomposition}
    There exists a subadditive nondecreasing sublinear function $d_0\in O(u)$ such that for every fixed $i,j\in\{1,\dots,n\}$, either all $d_0$-separated $i$-horizontal pairs are $j$-compressed or all such pairs are $j$-uncompressed. 
\end{lemma}

\begin{remark}
    Parts of the proof below (specifically in step 2) could have been made slightly simpler under the assumption that the $X_i$ and $Y_j$ are uniformly discrete. While natural in our context, we do not impose this assumption here.
\end{remark}

\begin{proof}
    To set notation, we use superscript for the factor index (i.e.\ $x^i:=\pi_i(x)$) and subscript to denote a sequence of points, for example $(x_k)_k$. The proof often focuses on one factor, say $X_j$, and we use $\widehat{X}_j:=\prod_{i\ne j}^n X_i$ to denote the remaining factors. Usually we only write the proof for the case $i=1$, in which case we simply write $\widehat{X}$. In these situations we sometimes write $x=(x^1,\widehat{x})$. For $\hat x \in \widehat X$, by $\hat x$-fiber we mean the set $X_1 \times \{\hat x \}$ and by $\widehat X$-fiber we mean some $\hat x$ fiber. 
    
    Fix $i,j\in\{1,\dots,n\}$.
    
    \paragraph{\textbf{Step $1$.}}  Let $x\in X$. We claim that there is a constant $d_0(x)$ such that for all $d\geq d_0(x)$, all $d$-separated $i$-horizontal pairs in the ball $B:=B\big(x,99d\big)$ are simultaneously $j$-compressed or simultaneously $j$-uncompressed. 
     
    First we show there is $d_0(x)$ such that for $d\geq d_0(x)$ there are no $d$-separated $i$-horizontal pairs in $B(x,99d)$ that are $j$-semi-compressed. Assume towards contradiction there were such pairs for all $d>0$. 
    These pairs give rise to a sequence of radii $d_k\rightarrow\infty$ and points $x_k,z_k\in B(x,99d_k)$ that are $i$-horizontal, $d(x_k,z_k)=d_k$, with $(x_k,z_k)$ being $j$-semi-compressed. Consider the asymptotic cone $X_\omega:=\operatorname{Cone}_\omega(X,x,d_k)$. It is $u$-admissible since the base point $x$ is fixed. Moreover, $x_k,z_k\in B(x,99d_k)$ so we have $x_\omega:=[(x_k)_k],z_\omega:=[(z_k)_k]\in X_\omega$. Since the pairs are $i$-horizontal, Lemma~\ref{lem: cone contradiction} yields a contradiction: $\omega$-almost surely,  $(x_k,z_k)$ is $l$-compressed for $l\ne i$ and $i$-uncompressed. We conclude that there is $d_0(x)$ such that $d>d_0(x)$ implies that inside $B(x,99d)$ there are no $i$-horizontal $d$-separated $j$-semi-compressed pairs. 
    
With a similar argument, we now show that there is $d_0(x)$ such that for every $d\geq d_0(x)$, every two pairs of $d$-separated $i$-horizontal points inside $B(x,99d)$ are $j$-compatible. Assume towards contradiction there is no such $d_0(x)$. One obtains a sequence of radii $d_k\to\infty$ and sequences $x_{1,k},z_{1,k},x_{2,k},z_{2,k}\in B(x,99d_k)$ such that for all $k$ and $\alpha\in\{1,2\}$, the pair $(x_{\alpha,k},z_{\alpha,k})$ is $i$-horizontal with $d_X(x_{\alpha,k},z_{\alpha,k})\in [d_k,99d_k)$ and $(x_{1,k},z_{1,k})$ is $j$-compressed while $(x_{2,k},z_{2,k})$ is $j$-uncompressed. We consider the $u$-admissible cone $X_\omega:=\operatorname{Cone}_\omega(X,x,d_k)$. Since $d_X(x_{\alpha,k},z_{\alpha,k})\in [d_k,99d_k)$, we see that $ x_{\omega,\alpha}:=[(x_{\alpha,k})_k],z_{\omega,\alpha}:=[(z_{\alpha,k})_k]\in X_\omega$. Since all pairs are $i$-horizontal, we conclude from Lemma~\ref{lem: cone contradiction} that if $l\ne i$, $\omega$-almost surely the pair $(x_{2,k},z_{2,k})$ is $l$-compressed. The assumption on $(x_{2,k},z_{2,k})$ forces $j=i$. On the other hand, Lemma~\ref{lem: cone contradiction} gives that  $\omega$-almost surely $(x_{1,k},z_{1,k})$ is $i$-uncompressed, so $j\ne i$ - a contradiction. 

We conclude the existence of a function $d_0$ as declared in the beginning of step 1: for all $d\geq d_0(x)$, all $d$-separated $i$-horizontal pairs in the ball $B:=B\big(x,99d\big)$ are simultaneously $j$-compressed or simultaneously $j$-uncompressed.

\paragraph{\textbf{Step~$2$}}    
    We now show that we can choose $d_0\in O(u)$. We start by taking, for each $x\in X$, the infimal $\tilde d_0(x)$ such that for all $d\geq \tilde d_0(x)$, all $i$-horizontal $d$-separated pairs in $B(x,99d)$ are $j$-compatible. Next define the radial nondecreasing function $d_0(|x|):=\sup_{y:|y|\leq|x|}\tilde d_0(y)$.
    
    A-priori $d_0:\mathbf{R}_{\geq 0}\rightarrow \mathbf{R}_{\geq 0} \cup \{\infty\}$, and we need to show it is in fact real valued (this is trivial if the spaces are uniformly discrete). Assume towards contradiction that $d_0$ is not bounded in some ball $B:=B(x_0,R)$, and get points $x_k\in B$ and real numbers $r_k\rightarrow\infty$ such that the balls $B(x_k,99r_k)$ contain two $i$-horizontal $r_k$-separated pairs $(x_{1,k},z_{1,k}), (x_{2,k},z_{2,k})$ that are $j$-incompatible - say that $(x_{1,k},z_{1,k})$ are $j$-compressed while $(x_{2,k},z_{2,k})$ are $j$-uncompressed. We consider the cone $X_\omega:=(X,x_0,r_k)$, that is $u$-admissible. Since $x_k\in B$ we have $d(x_k,x_0)\leq R$, so $x_{\omega,\alpha}:=[(x_{\alpha,k})],z_{\omega,\alpha}:=[(z_{\alpha,k})]\in X_\omega$. By Lemma~\ref{lem: cone contradiction} we see that $j$ cannot equal $i$ because of the assumption on $(x_{1,k},z_{1,k})$, and $j$ must equal $i$ because of the assumption on $(x_{2,k},z_{2,k})$ - a contradiction.

    The function $d_0$ is therefore a radial, real valued, bounded near $0$, and nondecreasing. We show it is $O(u)$.  Assume towards contradiction that there is a sequence $x_k\rightarrow \infty$ so that $\lim_\omega\frac{u(x_k)}{d_0(x_k)} = 0$. By the construction of $d_0$ we may take $x_k'$ such that $\vert x_k'\vert\leq \vert x_k\vert$, $d_0(x_k')=d_0(x_k)$, and for which $B\big(x_k',99d_0(x_k')\big)$ contains two $\big(d_0(x_k')-1\big)$-separated $i$-horizontal pairs $(x_{1,k},z_{1,k}), (x_{2,k},z_{2,k})$ that are $j$-incompatible. Since $u$ is nondecreasing, the assumption $\lim_\omega \frac{u(x_k)}{d_0(x_k)} = 0$ implies that the cone with moving base points $X_\omega:=\operatorname{Cone}_\omega\big(X,x_k',d_0(x_k')\big)$ is a $u$-admissible asymptotic cone. Since $x_{\alpha,k},z_{\alpha,k}\in B\big(x_k',99d_0(x_k)\big)$, we see that $x_{\omega,\alpha}:=[(x_{\alpha,k})],z_{\omega,\alpha}:=[(z_{\alpha,k})]\in X_\omega$. Once again, Lemma~\ref{lem: cone contradiction} yields a contradiction.

    \paragraph{\textbf{Step~$3$}}

We have established the following: there is $d_0\in O(u)$ such that for every $x\in X$ and every $d\geq d_0(x)$, the property of whether a $d$-separated $i$-horizontal pair inside $B=B(x,99d)$ is $j$-uncompressed or $j$-compressed depends only on $x$, $i$, $j$ and $d$, and not on the specific pair (in particular such a pair cannot be $j$-semi-compressed). We call it the $j$-type of $B$. Since $d_0\in O(u)$ and nondecreasing we may assume that $d_0$ is also subadditive. We can right away conclude that there are no $d_0$-separated $i$-horizontal $j$-semi-compressed pairs at all in $X$: such a pair $(x,x')$ must admit $d(x,x')\geq d_0(x)$
and therefore for $d=d(x,x')\geq d_0(x)$, $x,x'$ are $d$-separated and $x'\in B(x,99d)$. Step~$1$ shows that $(x,x')$ is not $j$-semi-compressed.

Our next goal is to remove any restriction on the radius of the balls $B(x,99d)$. Fix $x\in X$. First observe that the $j$-type of $d_0$-separated $i$-horizontal pairs which involve $x$ depend only on $x$, as they all have the same type as the type of the ball $B:=B\big(x,99d_0(x)\big)$. Explicitly, let $(x,x')$ be a $d_0$-separated $i$-horizontal pair. If $x'\in B$, $(x,x')$ has the same $j$-type as $B$. If $x\notin B$, we may connect $x$ and $x'$ by a geodesic $\gamma:=[x,x']$, and subdivide this geodesic into subintervals for which the endpoint have the same type as that of $B$. Precisely, the pair $\Big(x,\gamma\big(98d_0(x)\big)\Big)$ is $d>d_0(x)$-separated and $i$-horizontal, so it has the $j$-type of $B$. Therefore for $R_1:=98d_0(x)$, the ball $B_1:=B(x,99 R_1)$ has the $j$-type of $B$. We can continue to enlarge the radius of the balls until we find one which contains $x'$, so indeed $(x,x')$ has the $j$-type of $B$, which is what we were set out to prove

The next challenge is to consider $i$-horizontal pairs that do not share a point. Let $(x,x')$ and $(z,z')$ be two $d_0$-separated $i$-horizontal pairs. Assume first that $x,x',z,z'\in X$ are pairwise $i$-horizontal, so all points are in the same fiber. This means that the $\widehat{X}_i$ component of $x,x',z,z'$ is equal to some fixed $\hat{x}\in \widehat{X}_i$. Since each $X_i$ can be assumed to have infinite diameter, we can find an auxiliary  point $p=(p_1,\hat{x})$ which is far enough so that $(p,w)$ is $d_0$-separated for all $w\in \{x,x',z,z'\}$. It is clear that the $j$-type of all points must be equal to that of $B\big(p,99d_0(p)\big)$. Therefore for every choice of $\hat{x}\in\widehat{X}_i$ we have a well defined $(i,j,\hat{x})$-type of an $i$-horizontal $d_0$-separated pair inside the fiber $X_1\times\{\hat{x}\}$.

To finish the proof we need to show that the above type does not depend on $\hat{x}$. Fix $\hat{x},\hat{z}\in \widehat{X}$, $\hat{x}\ne \hat{z}$.
Since $X_i$ is of infinite diameter, there is a sequence $m_k$ tending to $\infty$ and point $w_k\in X_i$ such that $d_{X_i}(0_i,w_k)=m_k$. Consider the points $x_k=\big(w_k,\widehat{x}\big), z_k=\big(w_k,\widehat{z}\big)\in X$, the points whose $X_i$-coordinate is $w_k$, and which lie in the $\hat{x}$ (for $x_k$) or $\hat{z}$ (for $z_k$) fiber. For all large enough $k$, $(x_0,x_k)$ and $(z_0,z_k)$ are $d_0$-separated $i$-horizontal pairs. The pairs $(x_0,x_k)$ are all in the same $X_i$ fiber, determined by $\hat{x}$. For all large enough $k\in\mathbb{N}$, they are $d_0$-separated, so for all choices of large enough $k,l$, the two $i$-horizontal $d_0$-separated pairs $(x_0,x_k)$ and $(x_0,x_l)$ are $j$-compatible -- they have the $(i,j,\hat{x})$-type. Symmetrically, $(z_0,z_k)$ and $(z_0,z_l)$ are $j$-compatible with the $(i,j,\hat{z})$-type. Consider the cone $X_\omega:=(X,x_0,m_k)$, which is $u$-admissible and contains $x_\omega:=[(x_k)_k],z_\omega:=[(z_k)_k]$. Clearly $x_\omega=z_\omega$ and, for the base point $x_\omega^0$ we have $d(x_\omega^0,x_\omega)=1$. Using Lemma~\ref{lem: cone contradiction} with $x_\omega^0$ as one of the points in the cone and $x_\omega=z_\omega$ as the other point, we conclude that for all large enough $k,l\in \mathbb{N}$, the pairs $(x_0,x_k)$ and $(z_0,z_l)$ are $j$-compatible. Hence  the $(i,j,\hat{x})$-type equals the $(i,j,\hat{z})$-type, proving the lemma. 
\end{proof}

For clarity of presentation we make the following notations.

\begin{definition}
    For $a,b\in\mathbf{R}$, denote $a\lesssim b$ if there is $M>0$ such that $a\leq Mb$, and $a\sim b$ if both $a\lesssim b$ and $b\lesssim a$ hold. Let  $(a_k)_k,(b_k)_k$ be sequences of real numbers. We denote
    \begin{enumerate}
         \item  $a_k\lesssim b_k$ if there is $M>0$ independent of $k$ such that $a_k\leq Mb_k$ for all large enough $k$, and $a_k\sim b_k$ if both $a_k\lesssim b_k$ and $b_k\lesssim a_k$ hold. 
        \item $a_k\lesssim_\omega b_k$ if there is $M>0$ independent of $k$ such that $a_k\leq Mb_k$ $\omega$-almost-surely, and $a_k\sim_\omega b_k$ if both $a_k\lesssim_\omega b_k$ and $b_k\lesssim_\omega a_k$ hold. 
    \end{enumerate}
\end{definition}

We will repeatedly and implicitly use the following easy facts for an $O(v)$-bilipschitz embedding:
\begin{enumerate}
    \item If $(z_k)$ is a sequence in $Z$ with $|z_k|\rightarrow \infty$, then $|z_k|\sim |f(z_k)|$.
    \item If $(w_k)$ is a sequence in $W$ with $|w_k|\rightarrow \infty$, and $|z_k|$ are such that $d_W\big(f(z_k),w_k\big)\leq v(w_k)$, then $|z_k|\sim |w_k|$. Such $z_k$ exist if $f$ is an $O(v)$-equivalence. 
\end{enumerate}

\begin{proof}[Proof of Theorem~\ref{thm:KKLgeneralizedNEW}]
   From Lemma~\ref{lem: simultaneously-compressed-decomposition} we conclude that for every $i$ there is a unique $j$ for which the projection on $Y_j$ is uncompressed for all $d_0$-separated $i$-horizontal pairs. To ease  notation we assume $j=i$ is that unique index, and fix from now on $i=j=1$. The proof is identical for all other $i$.

\textbf{The $O(u)$-embeddings.}   We denote  $\widehat{X}:=\prod_{l=2}^n X_l$, and write points $x\in X$ as a pair $(z,\widehat{x})\in X_1\times \widehat{X}$. Fixing $\widehat{x}\in \widehat{X}$ we obtain a map $\phi_{\widehat{x}}:X_1\rightarrow Y_1$ given by $z\mapsto \pi_{Y_1}\big(\phi(z,\widehat{x})\big)$. We prove that these are all $(L_1,u_{\widehat{x}})$-bilipschitz embeddings, for $u_{\widehat{x}}=u+d_0+(u+d_0)(|\widehat{x}|_{\widehat{X}})\in O(u)$, where $d_0\in O(u)$ is given by Lemma~\ref{lem: simultaneously-compressed-decomposition}. Indeed, let $z,z'\in X_1$ be a two points such that the pair $(z,\widehat{x}),(z',\widehat{x})$ is $d_0$-separated. By Lemma~\ref{lem: simultaneously-compressed-decomposition} this pair is $1$-uncompressed and 
    \begin{equation*}
        L_1^{-1}d(z,z')\leqslant d\Big(\pi_{Y_1}\big(\phi(z,\widehat{x})\big),\pi_{Y_1}\big(\phi(z',\widehat{x})\big)\Big)\leqslant L_1d(z,z').
    \end{equation*}
    If the pair $(z,\widehat{x}),(z',\widehat{x})$ is not $d_0$-separated, then subadditivity of $d_0$ gives $d(z,z')=d_{X_1}(z,z')\leqslant d_0(|z|\vee |z'|)+d_0(|\widehat{x}|)$. Recalling that $L_1^{-1}<1$, this gives a lower bound:
    \begin{equation*}
       L_1^{-1}d(z,z')-d_0(|z|\vee |z'|)-d_0(|\widehat{x}|)\leqslant 0\leqslant d\Big(\pi_{Y_1}\big(\phi(z,\widehat{x})\big),\pi_{Y_1}\big(\phi(z',\widehat{x})\big)\Big). 
    \end{equation*}
    The Pythagorean formula yields the upper bound:
    \begin{equation*}
       d\Big(\pi_{Y_1}\big(\phi(z,\widehat{x})\big),\pi_{Y_1}\big(\phi(z',\widehat{x})\big)\Big)\leq d\big(\phi(z,\widehat{x}),\phi(z',\widehat{x})\big)\leq Ld(z,z')+u(|z|\vee |z'|)+u(|\widehat{x}|).
    \end{equation*}
    We conclude that $\phi_{\widehat{x}}:X_1\rightarrow Y_1$ is an $(L,u_{\widehat{x}})$-bilipschitz embedding.

\textbf{$O(u)$-surjectivity.} 
We claim that there is $v\in O(u)£$ such that for all $\hat{x}\in \widehat{X},y^1\in Y_1$, it holds that $d_{Y_1}(\operatorname{Im}\phi_{\hat{x}},y)\leq v(|\hat{x}|\vee|y^1|)$. Assume towards contradiction that this is not the case. This means there are sequences $\widehat{x_k}\in \widehat{X},y_k^1\in Y_1$ for which $\lim\frac{u(|\widehat{x_k}|\vee|y_k^1|)}{d_{Y_1}(\operatorname{Im}\phi_{\hat{x}},y_k^1)}=0$.

We make the following construction and notation:
\begin{itemize}
    \item Let $x_k^1\in X_1$ be points realizing $d_{Y_1}\big(\phi_{\widehat{x_k}}(x_k^1),y_k^1\big)\leq d_{Y_1}(\operatorname{Im}\phi_{\widehat{x_k}},y_k^1)+1$.
    \item $\sigma_k:=d_{Y_1}\big(\phi_{\widehat{x_k}}(x_k^1),y_k^1\big)$. 
    \item $x_k:=(x_k^1,\widehat{x_k})$.
    \item $\widehat{y_k}:=\pi_{\widehat{Y}}\circ\phi(\widehat{x_k},x_k^1)$. 
    \item $y_k:=(y_k^1,\widehat{y_k})$.
    \item Let $z_k=(z_k^1,\widehat{z_k})$ be points realizing $d_Y\big(\phi(z_k),y_k\big)\leq u(|y_k|)$.
\end{itemize}
Recall that $\phi_{\widehat{x_k}}$ is an $(L,u_{\widehat{x_k}})$-bilipschitz embedding, where $u_{\widehat{x_k}}=u+d_0+(u+d_0)(\widehat{x_k})$ and $d_0\in O(u)$. To ease notation, we will write $u+d_0=g\in O(u)$. We have
$$\frac{1}{L}|x_k^1|-g(|x_k^1|)-g(\vert \widehat{x_k} \vert)\leq d_{Y_1}\big(\phi_{\widehat{x_k}}(x_k^1), \phi_{\widehat{x_k}}(0_{X_1})\big)\leq d_{Y_1}\big(\phi_{\widehat{x_k}}(x_k^1), 0_{Y_1})\big)+d_{Y_1}\big(0_{Y_1},\phi_{\widehat{x_k}}(0_{X_1})\big).$$
By construction of $x_k$, the first term is bounded above by $\sigma_k+|y_k^1|$. For the second term we have $d_{Y_1}\big(0_{Y_1},\phi_{\widehat{x_k}}(0_{X_1})\big)=d_{Y_1}\big(\pi_{Y_1}\circ \phi(0_X),\pi_{Y_1}\circ\phi(0_{X_1},\widehat{x_k})\big)\lesssim |\widehat{x_k}|$. We conclude: 
\begin{equation}\label{eq: bounding size of _k^1}
    |x_k^1|\lesssim \sigma_k+|y_k^1|+|\widehat{x_k}|.
\end{equation}
We also have the following estimates:
\begin{equation}\label{eq: size estimates in surjectivity construction}
\begin{split}
      |\widehat{y_k}|\lesssim &|x_k^1|.\\
      |z_k|\sim &|(y_k^1,\widehat{y_k})|\lesssim |x_k^1|+|y_k^1|\lesssim |y_k^1|+\sigma_k+|\widehat{x_k}|.
\end{split}
\end{equation}
By our assumption that $\lim_k\frac{u(|\widehat{x_k}|\vee|y_k^1|)}{\sigma_k}=0$, we therefore have: 
\begin{equation}\label{eq: size estimates of u in surjectivity construction}
      \lim_k\frac{u(|x_k|)+u(|z_k|)+u(|y_k|)}{\sigma_k}=0.
\end{equation}
Moreover, using the defining inequality of $\phi$ and the formula for the distance in product spaces, we have:  
\begin{equation}\label{eq: bounding x_k-z_k}
\begin{split}
    \frac{d_X(x_k, z_k)}{L}-u(\vert z_k \vert \vee \vert x_k \vert)&\leq d_Y\big(\phi(x_k),\phi(z_k)\big)\\
   & \lesssim d_{Y_1}\big(\phi_{\widehat{0}}(x_k^1),\pi_{Y_1}\circ\phi(z_k)\big)+d_{\widehat{Y}}\big(\pi_{\widehat{Y}}\circ \phi(x_k),\pi_{\widehat{Y}}\circ \phi(z_k)\big)\\
   & \leq d_{Y_1}\big(\phi_{\widehat{0}}(x_k^1),y_k^1\big)+d_{Y_1}\big(y_k^1,\pi_{Y_1}\circ\phi(z_k)\big)+d_{\widehat{Y}}\big(\widehat{y_k},\pi_{\widehat{Y}}\circ \phi(z_k)\big)\\
   & \leq \sigma_k+u(y_k)+u(y_k)\lesssim \sigma_k,
\end{split}
\end{equation}
where the last inequality is a consequence of our assumption that $\lim \frac{u(|y_k^1|)}{\sigma_k}=0$.

Combining \eqref{eq: bounding x_k-z_k} with the estimates in~\eqref{eq: size estimates in surjectivity construction} we conclude:
\begin{equation}\label{eq: z_k^1,0 and z_k are in the cone}
    d_X(x_k, z_k)\lesssim \sigma_k+g(|z_k|\vee|x_k|)\lesssim \sigma_k+g(|y_k^1|)+g(|\widehat{x_k}|)+g(\sigma_k)\lesssim \sigma_k,
\end{equation}
where again the last inequality follows from \eqref{eq: size estimates of u in surjectivity construction} and $g\in O(u)$. Consider the cone $X_\omega:=\operatorname{Cone}(X,x_k,\sigma_k)$. Inequality~\eqref{eq: size estimates of u in surjectivity construction} establishes that $X_\omega$ is $u$-admissible. Moreover, \eqref{eq: z_k^1,0 and z_k are in the cone} ensures that the points $z_\omega:=[(z_k)_k]$ and $z_\omega^x:=[(z_k^1,\widehat{x_k})_k]$ both lie in $X_\omega$. Therefore by Lemma~\ref{lem: cone contradiction} we get 
$$\lim_\omega \frac{1}{\sigma_k}d_{Y_1}\big(\pi_{Y_1}\circ\phi(z_k^1,\widehat{x_k}),\pi_{Y_1}\circ\phi(z_k^1,\widehat{z_k})\big)=0,$$
and therefore 
$$\lim_\omega\frac{1}{\sigma_k}d_{Y_1}\big(\phi_{\widehat{x_k}}(z_k^1),y_k^1\big)=\lim_\omega\frac{1}{\sigma_k}d_{Y_1}\big(\pi_{Y_1}\circ\phi(z_k),y_k^1\big)\leq \lim\frac{u(|y_k|)}{\sigma_k}=0$$
(again, the last step is by Inequality~\eqref{eq: size estimates of u in surjectivity construction}). We reach a contradiction to the definition of $\sigma_k<d(\operatorname{Im}\phi_{\widehat{x_k}},y_k^1)+1$. We conclude the existence of $v\in O(u)$ as claimed.

\textbf{Sublinear Control.} We show there are sublinear functions $h,g$ with $g\in O(u)$ such that for $\widehat{x},\widehat{w}\in \widehat{X}$ and $x^1\in X_1$, we have $d_{Y_1}\big(\phi_{\widehat{x}}(x^1),\phi_{\widehat{w}}(x^1)\big)\leq h\big(\vert\widehat{x}\vert\vee\vert\widehat{w}\vert)\vee g(\vert x^1\vert)$. 
We show it for $\widehat{w}=\widehat{0}$, which by the triangle inequality implies the general case. 

Assume towards contradiction that there are sequences $(x_k^1)_k$ in $X_1$ and $(\widehat{x_k})_k$ in $\widehat{X}$ such that $\lim_\omega\frac{u(\vert x_k^1\vert)}{\sigma_k}=0$ and $\lim_\omega\frac{\sigma_k}{\vert \widehat{x_k} \vert}>0$, where  $\sigma_k=d_{Y_1}\big(\phi_{\widehat{x_k}}(x_k^1),\phi_{\widehat{0}}(x_k^1)\big)$. Denote $x_k=(x_k^1,\widehat{x_k})$ and consider the cone $X_\omega=\operatorname{Cone}(X,x_k,\sigma_k)$. The assumption on $\sigma_k$ once again implies that $X_\omega$ is $u$-admissible, and that $x_\omega^0:=[(x_k^1,\widehat{0})_k]\in X_\omega$. Lemma~\ref{lem: cone contradiction} implies that $\lim_\omega \frac{\sigma_k}{\sigma_k}=0$, a contradiction. We conclude the existence of $g\in O(u)$ and a sublinear $h$ such that $d_{Y_1}\big(\phi_{\widehat{x}}(x^1),\phi_{\widehat{w}}(x^1)\big)\leq h\big(\vert\widehat{x}\vert\vee\vert\widehat{w}\vert)\vee g(\vert x^1\vert)$ for all $\widehat x$ and $\widehat w$ in $\widehat X$, as claimed. Moreover, we have 
$$d_{Y_1}\big(\phi_{\widehat{x}}(x^1),\phi_{\widehat{w}}(x^1)\big)\leq d_Y\big(\phi(x^1,\widehat{x}),\phi(x^1,\widehat{w})\big)\leq Ld_{\widehat X}(\widehat x, \widehat w)+u(\vert \widehat{x} \vert \vee \vert \widehat{w}\vert) + u(\vert x^1\vert),$$
hence we can take $h$ to be bounded near $\widehat{0}$. By defining $\tilde{h}(r):\sup_{x\leq r}=h(r)$ we may assume $h$ is nondecreasing, and by Lemma~\ref{lem:dominating-sublinear-by-subadditive} we conclude that one can take $h$ to be subadditive and nondecreasing. We let $v=h+g$, finishing the proof. \qedhere

\end{proof}

\subsection{Proof of Theorem~\ref{thm:KKLgeneralized}} \label{sec: proving theorem KKL}

We recall the three definitions below from~\cite{KKL}.

\begin{definition}
Let $X$ be a geodesic metric space. We say that $p,q\in X$ are in the same \emph{leaf} if there is a continuous path $\gamma:I\rightarrow X$ joining $p$ and $q$ such that every other continuous path which joins $p$ and $q$ contains $\gamma$. Being in the same leaf is an equivalence relation, and every leaf of $X$ is a closed convex subset. The space $X$ is called \emph{type} I if all its leaves are geodesically complete trees which branch everywhere. 
\end{definition}

\begin{definition}
    A geodesic metric space is of \emph{type} II if it is a thick irreducible Euclidean building with transitive affine Weyl group of rank $r\geq 2$.
\end{definition}

\begin{definition}
    A space $X$ is of \emph{coarse type} I (resp. II) if it is geodesic and every asymptotic cone of $X$ is of type I (resp. II).
    \end{definition}

Symmetric spaces of noncompact type and high rank are coarse type II \cite{KleinerLeebQI}, and so are Euclidean buildings of high rank and cocompact affine Weyl group. Examples of spaces of coarse type I are Gromov hyperbolic groups whose Gromov boundary contains at least 3 points. These examples are enough for all our purposes in this paper, which arise in Section~\ref{sec: completely solvable groups and corollaries} in the proofs of Corollary~\ref{cor:KKL-symmetric} and Theorem~\ref{th:technical-KKL}.

In order to include Euclidean factors (or, more generally, factors with asymptotic cones homeomorphic to Euclidean space) we will need a slight variation of Theorem~\ref{thm:KKLgeneralizedNEW}. 

\begin{proposition}\label{prop: product with Euclidean factors}
    Let $X:=\overline{X}\times Z,\ Y:=\overline{Y}\times W$ be two products of geodesic metric spaces, and let $\phi:X\rightarrow Y$ be an $(L,u)$-bilipschitz equivalence. Assume that in any pair of $u$-admissible cones $X_\omega,Y_\omega$, there is a homeomorphism $\overline{\Phi}:\overline{X}_\omega\rightarrow \overline{Y}_\omega$ such that for the induced cone map $\Phi:=\operatorname{Cone}(\phi):X_\omega\rightarrow Y_\omega$, the following diagram commutes: 
     \[
    \begin{tikzcd}
X_\omega \ar[r, "\Phi"] \ar[d,"\overline{\pi}_\omega"] & Y_\omega \ar[d,"\overline{\pi}_\omega"] \\
\overline{X}_\omega \ar["\overline{\Phi}",r] & \overline{Y}_\omega
\end{tikzcd}
\]
Then there is an $O(u)$-bilipschitz equivalence $\overline{\phi}:\overline{X}\rightarrow \overline{Y}$ such that the diagram 
\[
    \begin{tikzcd}
X \ar[r, "\phi"] \ar[d,"\overline{\pi}"] & Y \ar[d,"\overline{\pi}"] \\
\overline{X} \ar["\overline{\phi}",r] & \overline{Y}
\end{tikzcd}
\]
commutes up to sublinear error $v$ that is nondecreasing and subadditive.
\end{proposition}

We omit the proof, which is very similar to the proof of Theorem~\ref{thm:KKLgeneralizedNEW}. 

\begin{remark}\label{rmk: omitted proof of Euclidean factors proposition}
    There is one essential difference between the proof of Proposition~\ref{prop: product with Euclidean factors} and that of Theorem~\ref{thm:KKLgeneralizedNEW}: in the former, the constants $L_1$ and $\epsilon$ in the definition of compressed and uncompressed pairs are not chosen freely. Rather, we can only guarantee that there are such constants $L_1>L$ and $\epsilon<L_1^{-1}$ for which the proof carries. The reason is that in the asymptotic cones, one does not have good enough control over the $Z$ factor. The result is that the multiplicative (Lipschitz) constant of $\overline{\phi}$ may be larger than $L$. We stress that $\overline{\phi}$ is still an $O(u)$-equivalence. 

\end{remark}

\begin{proof}[Proof of Theorem~\ref{thm:KKLgeneralized}]
We denote $\overline{X}:=\prod_{i=1}^n X_i, \overline{Y}:=\prod_{j=1}^m Y_j$. The map $\phi:X\rightarrow Y$ is an $O(u)$-equivalence, therefore in any $u$-admissible cone $X_\omega$, the map $\phi_\omega$ is an $L$-bilipschitz homeomorphism (Lemma~\ref{lem:go-through-cones}).  The assumption on the factors then allows us to use Theorem~\ref{thm: topological splitting}, and conclude that $p=q,n=m,$ and that the map $\overline{\phi}_\omega: \pi_{\overline{Y}_\omega}\circ\phi_{\omega{\restriction \{0_{\mathbf{R}^p}\}\times \overline{X}_\omega}}:\overline{X}_\omega\rightarrow \overline{Y}_\omega$ is a homeomorphism that preserves the product structure. 

We conclude that for $X,Y$ and $\phi$, the hypotheses of Proposition~\ref{prop: product with Euclidean factors} are met, and we obtain an $O(u)$-equivalence $\overline{\phi}:\overline{X}\rightarrow \overline{Y}$ such that 
\[
    \begin{tikzcd}
X \ar[r, "\phi"] \ar[d,"\overline{\pi}"] & Y \ar[d,"\overline{\pi}"] \\
\overline{X} \ar["\overline{\phi}",r] & \overline{Y}
\end{tikzcd}
\]

commutes up to a sublinear error $\overline v$ that is subadditive and nondecreasing. The triple consisting of $\overline{x}, \overline{Y}$ and $\overline{\phi}:\overline{X}\rightarrow \overline{Y}$ admit the hypotheses of Theorem~\ref{thm: topological splitting}, hence in any $u$-admissible cones $\overline{X}_\omega,\overline{Y_\omega}$, the map $\overline{\phi}_\omega:\overline{X}_\omega\rightarrow\overline{Y}_\omega$ preserves the product structure. Therefore the hypotheses of Theorem~\ref{thm:KKLgeneralizedNEW} are met, and we obtain for all $i \in \{1, \ldots, n \}$, $O(u)$-equivalences $\overline{\phi}_i:X_i\rightarrow Y_i$ such that the diagrams 
 \[ \begin{tikzcd}
 \overline{X}=\prod_{i=1}^n X_i \ar[r, "\phi"] \ar[d,"\pi_i"] & \prod_{i=1}^n Y_i=\overline{Y} \ar[d,"\pi_{\sigma(i)}"] \\
X_i \ar["\phi_i",r] & Y_{\sigma (i)}
\end{tikzcd} \]
commute up to sublinear error $v_i$ that is subadditive and nondecreasing. Finally, the exterior square in the diagram 
 \[ \begin{tikzcd}
 X \ar[r, "\phi"] \ar[d, "\overline{\pi}"] & Y \ar[d, "\overline{\pi}"]\\
 \overline{X} \ar[r, "\overline{\phi}"] \ar[d,"\pi_i"] & \overline{Y} \ar[d,"\pi_i"] \\
X_i \ar["\overline{\phi}_i",r] & Y_i
\end{tikzcd} \]

commutes up to a sublinear error $v' = \overline v + \max_i v_i$, as can be checked as follows: for every $x \in X$,
\begin{align*}
    d(\pi_i \circ \overline \pi \circ \phi(x), \overline{\phi}_i \circ \pi_i \circ \overline \pi(x))  & \leqslant   d(\pi_i \circ \overline \pi \circ \phi(x), \pi_i \circ \overline \phi \circ \overline \pi (x)) \\
    & + d(\pi_i \circ \overline \phi \circ \overline \pi (x), \overline{\phi}_i \circ \pi_i \circ \overline \pi(x)) \\
    & \leqslant   d(\overline \pi \circ \phi(x), \overline \phi \circ \overline \pi (x)) + v_i(\overline \pi(x)) \\
    & \leqslant \overline v(\vert x \vert) + \max_{i \in \{1,\ldots, n\}} v_i(\vert x \vert).
\end{align*}
where we used the triangle inequality and then the fact that the projections $\overline \pi$ and $\pi_i$ are distance-decreasing.\footnote{The real point used here is that the projections $\overline \pi$ and $\overline \pi_i$ are $o(r)$-lipschitz and thus the diagram commutes in the $o(r)$-lipschitz category defined by Cornulier in  \cite{cornulier2017sublinear}.}
\end{proof}

\section{Completely solvable groups, Corollary~\ref{cor:KKL-symmetric} and Theorem~\ref{cor:KKLgeneralized}}\label{sec: completely solvable groups and corollaries}

\subsection{Some preliminaries on Cornulier's reduction}
\label{sec:cornulier-reduction-partial-converse}
As announced in the beginning of the introduction, the quasiisometry classification of connected Lie groups amounts to that of the completely solvable ones.
Given a simply connected solvable Lie group $G$ there are several possible, equivalent definitions for $\rho_0(G)$. We give below  that of Jablonski, building on Gordon and Wilson.

\begin{proposition}[{\cite[\S4.1 and \S4.2]{JabMax}, after \cite{GordonWilson}}]
    Let $G$ be a simply connected solvable Lie group. There exists a (possibly non-unique) left-invariant metric  $g_{\operatorname{max}}$ on $G$ whose isometry group contains a transitive completely solvable group $G_0$. Moreover, the group $G_0$ obtained in this way is unique up to isomorphism and it does not depend on $g_{\max}$.
\end{proposition}

\begin{definition}
    Let $G$ be a simply connected solvable Lie group. We define $\rho_0(G)$ as $G_0$. We say that a group is in the class $(\mathcal C_0)$ if $G = \rho_0(G)$, that is, if $G$ is completely solvable.
    \end{definition}

It is clear that the groups $G$ and $G_0$ are quasiisometric, being closed co-compact subgroups of the isometry group $\widehat G$ of $g_{\operatorname{max}}$. They are commable in the terminology of \cite{CornulierQIHLC}. The role of the group $\widehat G$ is played by the group denoted $G_3$ in Cornulier's treatment (\cite[Lemme 1.3]{Cornulieraspects}, summarizing \cite{CornulierDimCone}).

\begin{definition}
Let $G$ be a group in the class $(\mathcal C_0)$.
The exponential radical $\operatorname{R}_{\exp} G$ of $G$ is the smallest normal subgroup $N$ of $G$ such that $G/N$ is nilpotent.
\end{definition}
The exponential radical was named by Osin \cite{OsinExprad} as it is the subgroup of exponentially distorted elements in $G$ (together with $1$).
We call $\dim G/\operatorname{R}_{\exp} G$ the rank of $G$.
If $\widehat G$ is a real semisimple Lie group with trivial center, writing an Iwasawa decomposition $\widehat G = KAN$ and setting $G=AN$, we recover that the real rank of $\widehat G$ is the rank of $G$. More generally, the rank as defined here is still the dimension of one (or any) Cartan subgroup of $G$.

\begin{definition}
Let $G$ be a completely solvable Lie group with exponential radical $N$.
Say that $G$ is in $(\mathcal C_1)$ if the extension
$1 \to N \to G \to G/N \to 1$ splits and the action of $G/N$ on $N$ is $\mathbf R$-diagonalizable.
\end{definition}

\begin{definition}
    Let $G$ be a completely solvable group with $N = \operatorname{R}_{\exp} G$, and set $H = G/N$.
Decompose $\phi = \operatorname{ad}\colon \mathfrak g \to \operatorname{Der}(\mathfrak n)$ into
\[ \phi = \phi_\delta + \phi_\nu \]
where $\phi_\delta$ is $\mathbf R$-diagonalisable and $\phi_\nu$ is nilpotent \cite{BbkiCartanAlg}.
Note that $\phi_{\delta}$ is zero when restricted to $\mathfrak n$, so that it is well-defined on $\mathfrak h$.
Let $\rho_1(G)$ be $N \rtimes H$, where $\mathfrak h$ acts on $\mathfrak n$ through $\phi_\delta$. We also write $\rho_1(\mathfrak g)$ for $\operatorname{Lie}(\rho_1(G))$.
\end{definition}

Although in this paper the focus is on $\rho_1$, we recall below the definition of $\rho_\infty$, a further reduction that we mentioned in the introduction.

\begin{definition}[Cornulier, \cite{CornulierCones11}]
    Let $G$ be a completely solvable group, and let $H =G/\operatorname{R}_{\exp G}$. The Lie algebra of $H$ has a filtration by its the derived central series. Define 
    \[ \mathfrak h_\infty = \bigoplus_{i>0} C^i \mathfrak h/C^{i+1} \mathfrak h \]
    with the brackets induced from those of $\mathfrak h$. The action of $H$ on $\operatorname{R}_{\exp} G$, after being factored through $H/C^2H\simeq H_{\infty}/C^2H_\infty$, lifts a new action of $H_\infty$ on $\operatorname{R}_{\exp} G$; define $\rho_{\infty}(G)$ as the corresponding semidirect product $\operatorname{R}_{\exp G} \rtimes H_\infty$.
\end{definition}

We will require the following theorem of Cornulier: 

\begin{theorem}[{Cornulier, \cite{CornulierCones11}}]
\label{th:Cornulier-thm}
Let $G$ be a completely solvable group, and let $H = G/\operatorname{R}_{\exp} G$.
Then
\begin{enumerate}
\item $G$ and $\rho_1(G)$ are $O(\log)$-bilipschitz equivalent.
    \item $H$ is an $O(\log)$-Lipschitz retract of $G$, more precisely: \label{item:cornulier-retract} $\pi: G \to H$ is $O(\log)$-Lipschitz.

    \item $G$ and $\rho_\infty(G)$ are $O(u)$-bilipschitz equivalent, where the function $u$ depends on $G$.
\end{enumerate}

\end{theorem}

\subsection{Warm-up: the groups in Example \ref{exm:four-dim-KKL} are not quasiisometric.} \label{sec: example 1.3}

We prove below that the groups in Example \ref{exm:four-dim-KKL} are not quasiisometric. The proof is less involved than that of Theorem~\ref{cor:KKLgeneralized} but the main idea already intervenes, so we give it before.

\begin{proposition}
    Let $\alpha \in (0,1)$.
    The groups $G^0_{4,9}$ and $\mathbf R \times G_{3,5}^{\alpha}$ are not quasiisometric, for any $\alpha \neq 1$.
\end{proposition}

\begin{proof}
    We first check that $\mathbf R \times G_{3,3}$ is the group in the class $(\mathcal C_1)$ associated to $G_{4,9}^0$ by \cite{CornulierCones11}, and so there is a $O(\log)$-bilipschitz equivalence $\phi \colon G_{4,9}^0 \to \mathbf R \times G_{3,3}$. 
    The Lie algebra $\mathfrak g = \mathfrak g_{4,9}^0$ of $G_{4,9}^0$ is a semidirect product $\mathfrak{heis} \rtimes \delta$ where $\mathfrak{heis}$ has a basis $(e_1,e_2,e_3)$ with $[e_1,e_2] = e_3$ and $\delta e_i = e_i$ if $i=1,3$ or $0$ if $i=2$.
    The derived subalgebra is $\mathfrak u = \langle e_1, e_3 \rangle$ and $[\mathfrak g, \mathfrak u] = \mathfrak u$, so that $\mathfrak u$ is the Lie algebra of the exponential radical.
    The matrices of $\operatorname{ad}_{e_3}$ and $\operatorname{ad}_{e_4}$ in the basis $(e_1,e_2)$ are, respectively,
    \[ 
    \begin{pmatrix}
        0 & 1 \\ 0 & 0 
    \end{pmatrix}
    \quad \text{and} \quad 
    \begin{pmatrix}
        1 & 0 \\ 0 & 1
    \end{pmatrix}.
    \]
    From the definition of $\rho_1$ it follows that, using the same basis and changing the brackets, $\rho_1(\mathfrak g_{4,9})^0$ is the Lie algebra $\mathfrak u \rtimes \langle e_4 \rangle \times \mathbf R \simeq \mathfrak g_{3,3} \times \mathbf R$.
    
    Assume now towards contradiction that  $G^0_{4,9}$ and $\mathbf R \times G_{3,5}^{\alpha}$ are quasiisometric; then by Cornulier's Theorem \ref{th:Cornulier-thm} there is a $O(\log)$-bilispchitz equivalence $\psi \colon \mathbf R \times G_{3,3} \to \mathbf R \times G_{3,5}^{\alpha}$.
    By Theorem~\ref{thm:KKLgeneralized} there is a sublinear bilipschitz equivalence $\underline \psi \colon  G_{3,3} \to  G_{3,5}^{\alpha}$. However, by \cite{pallier2019conf}, $O(u)$-bilipschitz equivalences between Heintze groups preserve the conformal dimensions of their Gromov boundaries. 
    $G_{3,3}$ and $G_{3,5}^{\alpha}$ are Heintze groups, and the conformal dimensions are
    \[ \operatorname{Cdim} \partial_\infty G_{3,3} =2, \, \quad \text{while} \qquad \operatorname{Cdim} \partial_\infty G_{3,5}^{\alpha} =1+1/\alpha>2. \]
    So we reach a contradiction.
\end{proof}

\subsection{Proof of Corollary~\ref{cor:KKL-symmetric} and Theorem~\ref{cor:KKLgeneralized}}
\label{subsec:proof-corollaries}

We will achieve both corollaries by means a technical result, Theorem \ref{th:technical-KKL}.
Before that we need some preparation.

\begin{definition}
    A Heintze group of diagonal type is a completely solvable Lie group $M \rtimes \mathbf R$, where $M$ is simply connected nilpotent and $t \in \mathbf R$ acts by $\exp(tD)$, where $D \in \operatorname{Der}(\mathfrak m)$ is diagonalizable and has strictly positive eigenvalues.
\end{definition}

If $H$ is a Heintze group of diagonal type, then $M$ as above is its nilradical.

\begin{theorem}[{Heintze, \cite{Heintze}}]
    A group in $(\mathcal C_1)$ is a Heintze group of diagonal type if and only if it carries a left-invariant metric of strictly negative curvature.
\end{theorem}

Some examples of Heintze groups of diagonal type come from rank one symmetric spaces. More generally, we define Iwasawa subgroups as follows.

\begin{definition}
    Let $G$ be a simple Lie group with trivial center. Let $G=KAN$ be an Iwasawa decomposition of $G$. We call $AN$ an Iwasawa subgroup of $G$; it is unique up to conjugation in $G$ so that we may speak of ``the Iwasawa subgroup of $G$''. 
    If $G$ has real rank one, its Iwasawa subgroup is a Heintze group of diagonal type, see \cite[\S 9.3]{PansuCCqi}.
\end{definition}

Given a metric space $X$ (or a group) and an admissible function $u$, we denote by $\operatorname{SBE}(X)^{O(u)}$ the group of self $O(u)$-bilipschitz equivalences of $X$.

\begin{definition}
Let $H$ be a Heintze group of diagonal type.
    We say that $H$ satisfies the strong pointed sphere property if for every admissible $u$ the group $\operatorname{SBE}(H)^{O(u)}$ is not transitive on $\partial_\infty H$. 
\end{definition}

The Iwasawa subgroups of simple Lie groups of rank one do not have the strong pointed sphere property, since the group of isometries of a rank one symmetric space acts transitively on its Gromov boundary.
Among the Heintze groups with an abelian nilradical, the ones of the latter kind turn out to be the only exceptions:

\begin{lemma}[{\cite[Lemma 4.1]{loglie}}]
    Let $H$ be a completely solvable Heintze group with an abelian nilradical. If $H$ is not isomorphic to a maximal completely solvable subgroup of $\mathrm{SO}(n,1)$ for any $n\geqslant 2$, then $H$ has the strong pointed sphere property.
\end{lemma}

We will also need the following lemma:

\begin{lemma}[{Consequence of \cite[Theorem 1.4]{GrayevskyRigidity}}]
    \label{lem:SBE-classification-higher-rank}
    Let $G_{\operatorname{II}}$ and $G'_{\mathrm{II}}$ be two semisimple Lie groups with trivial centers, no factors of rank one and no compact factors. If $G_{\operatorname{II}}$ and $G'_{\mathrm{II}}$ are $O(u)$-bilipschitz equivalent, then they are isomorphic.
\end{lemma}

\begin{proof}
Let $X$ be the Riemannian symmetric space associated to $G_{\mathrm{II}}$.   It follows from \cite{GrayevskyRigidity} (with $X = X_0$) that the group $\operatorname{SBE}(G_{\mathrm{II}})^{O(u)}=\operatorname{SBE}(X)^{O(u)}$ is isomorphic to $G_{\mathrm{II}}$.
Similarly, $\operatorname{SBE}(G'_{\mathrm{II}})^{O(u)}$ is isomorphic to $G'_{\mathrm{II}}$.
Assuming that $G_{\operatorname{II}}$ and $G'_{\mathrm{II}}$ are $O(u)$-bilipschitz equivalent, the groups $\operatorname{SBE}(G_{\mathrm{II}})^{O(u)}$ and $\operatorname{SBE}(G'_{\mathrm{II}})^{O(u)}$ are isomorphic, concluding the proof.
\end{proof}

\begin{remark}
    The way Lemma \ref{lem:SBE-classification-higher-rank} is deduced from \cite{GrayevskyRigidity} is the same the way the  QI classification was deduced in \cite{KleinerLeebQI}. 
    One may nevertheless wish for a direct way of proving this.
\end{remark}

\begin{theorem}\label{th:technical-KKL}

    Let $m$, $n$, $m'$ and $n'$ be nonnegative integers.
    Let $G$, and $G'$ be two real semisimple Lie groups with trivial center and no compact factors. 
    Let $\{ H_i \}_{1 \leqslant i \leqslant m}$ and $\{H'_j\}_{1 \leqslant j \leqslant m'}$ be families of diagonal Heintze groups satisfying the strong pointed sphere property.
    Let $R$ and $R'$ be simply connected nilpotent Lie groups with dimensions $n$ and $n'$ respectively as real Lie groups.
    Write $G=KAN$, and $G' = K'A'N'$.
    Form the following solvable groups $P$ and $P'$ :
    \begin{equation*}
        P  = R \times AN \times \prod_{i=1}^m H_i  \; \text{and}\
        P'  = R' \times A'N' \times \prod_{j=1}^{m'} H'_j.
    \end{equation*}
    Let $S$ and $S'$ be completely solvable groups. Assume that $S$ and $S'$ are $O(u)$-bilipschitz equivalent for some subadditive function $u$, that $\rho_1(S) = P$ and that $\rho_1(S') = P'$.
    Then 
    $m=m'$, $n=n'$, $G$ and $G'$ are isomorphic, and there is a bijection $\sigma$ of $\{1, \ldots, m \}$ such that $H_i$ is $O(u \vee \log)$-equivalent to $H_{\sigma(i)}$ for all $i$.

\end{theorem}

\begin{remark}
    The strong pointed sphere property is known to hold for some Heintze groups with nonabelian nilradicals \cite[Remark 9]{loglie}. So the assumptions of Theorem \ref{th:technical-KKL} are indeed weaker than those of Theorem \ref{cor:KKLgeneralized} (see also Remark~\ref{rem:end-remark} for concrete examples). Moreover, instead of the Euclidean factors of Theorem~\ref{cor:KKLgeneralized}, Theorem~\ref{th:technical-KKL} allows more general $R$ and $R'$. 
\end{remark}

\begin{proof}[Proof of Corollary \ref{cor:KKL-symmetric} assuming Theorem \ref{th:technical-KKL}]
    Let $L$ and $L'$ be semisimple Lie group with no compact factors such that $X = L/K$ and $Y = L'/K'$. Set $S = AN$, $S' = A'N'$, $n=n'=m=m'=0$. By assumption, $S$ and $S'$ are sublinear bilipschitz equivalent. By Lemma~\ref{lem:dominating-sublinear-by-subadditive}, $S$ and $S'$ are $O(u)$-bilipschitz equivalent for some subadditive function $u$. Applying the theorem to the pair $(S,S')$, it follows that $L=G$ and $L'=G'$ are isomorphic, so that $X$ and $Y$ are pluriisometric.
\end{proof}

\begin{proof}[Proof of Theorem \ref{cor:KKLgeneralized} assuming Theorem \ref{th:technical-KKL}]
    Let $S$, $S'$, $P$ and $P'$ be as in the assumptions of the Theorem \ref{cor:KKLgeneralized}.
    Then $S$, $S'$, $P$ and $P'$ also satisfy the assumptions of Theorem \ref{th:technical-KKL}, since the groups $H_i$ and $H'_j$ satisfy the strong pointed sphere property as recalled above. 
    By Theorem \ref{th:technical-KKL}, the groups $\mathbf R^n \times AN$ and $\mathbf R^{n'} \times A'N'$ are isomorphic, while, after possibly reindexing the groups $H'_j$, $H_i$ is $O(u \vee \log)$-equivalent to $H'_i$ for all $i$. Now by the main theorem of \cite{pallier2019conf}, $H_i$ and $H'_i$ are isomorphic for all $i$, and so $P$ and $P'$ are isomorphic.
\end{proof}

We now come to the proof of Theorem~\ref{th:technical-KKL}. Applying Theorem \ref{thm:KKLgeneralized} to $P$ and $P'$ we will get $O(\log)$-equivalences between the direct factors of $P$ and $P'$. The proof consists in establishing which factors can be paired with one another. We prove that simple factors of $G$ and $G'$ must pair to one another, whilst preserving the $\mathbf{R}$-rank one vs. high $\mathbf{R}$-rank distinction; the strong pointed sphere property further allows to ensure that the remaining factors, namely the Heintze groups of diagonal type with this property, are not paired with the $\mathbf R$-rank one factors which do not have it. 

\begin{proof}[Proof of Theorem \ref{th:technical-KKL}]
    Let $n,m,n',m',S,S'$, $P, P'$, $R$, $R'$, $G$ and $G'$ be as in the assumptions of Theorem \ref{th:technical-KKL}. 
    Decomposing $G$ and $G'$ into products of simple factors, we find that there exists $p$, $p'$, $q$ and $q'$ so that 
    \begin{equation*}
        G = \prod_{i=1}^{p+q} G_i \quad \text{and} \quad
        G' = \prod_{j=1}^{p'+q'} G_j,
    \end{equation*}
    where $G_i$ has $\mathbf R$-rank one for $1 \leqslant i \leqslant p$ and $\mathbf R$-rank at least two for $p+1 \leqslant i \leqslant p+q$, similarly for $G'_j$.
    Let $\ell = p+m$ and $\ell' = p'+m'$. For $i \in \{m+1, \ell\}$ define $H_i:= A_{i-m} N_{i-m}$, where $G_k = K_k A_k N_k$ for $1 \leqslant k \leqslant p$; similarly, define $H'_j$ for $j$ in $\{m'+1, \ell' \} $. 
    In this way, $\{H_i\}_{i=1}^m$, resp.\ $\{H_j\}_{j=1}^{m'}$ are the Heintze groups that were originally present in the statement, and $\{H_i\}_{i=m+1}^{m+p}$, resp.\ $\{H_j\}_{j=m'+1}^{m'+p'}$ are the Iwasawa subgroups of the rank $1$ factors of $G$, resp. $G'$ (which are also diagonal type Heintze groups).
    Especially, the groups $H_i$ for $1\leqslant i \leqslant m$ and $H'_j$ for $1\leqslant j \leqslant m'$ have coarse type I.

   It follows from Theorem~\ref{th:Cornulier-thm}(1) that the groups $S$ and $P$ on the one hand, $S'$ and $P'$ on the other hand, are $O(\log)$-bilipschitz equivalent.
    By assumption, $S$ and $S'$ are $O(u)$-bilipschitz equivalent, especially they are $O(u \vee \log)$-bilipschitz equivalent.
    So $P$ and $P'$ are $O(u \vee \log)$-bilipschitz equivalent.
    Note that the coarse type I factors of $P$, resp.\ $P'$ are the groups $H_i$, resp.\ the groups $H'_j$, while the coarse type II factors are of $P$, resp.\ $P'$ are the groups $A_kN_k$ for $k\in \{p+1,\ldots, q\}$, resp.\ $A'_sN'_s$ for $s \in \{p'+1,\ldots, p'+q'\}$, which are respectively quasiisometric to the simple Lie groups $G_k$ and $G'_s$ of high $\mathbf R$-rank.
    
    Applying Theorem \ref{thm:KKLgeneralized} to $P$ and $P'$ we find that $n=n'$, $q=q'$, $\ell = \ell'$, and there are two bijections $\sigma_{\mathrm{I}}$ of $\{1,\ldots, \ell\}$ and $\sigma_{\mathrm{II}}$ of $\{ p+1,\ldots p+q \}$ so that
    \begin{enumerate}[(I)]
        
        \item \label{item:coarse-typeI}  $H_{i}$ is $O(\log)$-bilipschitz equivalent to $H_{\sigma_{\mathrm{I}}(i)}$ for all $1 \leqslant i \leqslant \ell$; 
        \item \label{item:higher-rank} $G_{i}$ is $O(\log)$-bilipschitz equivalent to $G_{\sigma_{\mathrm{II}}(i)}$ for all $p<i\leqslant p+q$; 
    \end{enumerate}

    From \eqref{item:higher-rank} and Lemma~\ref{lem:SBE-classification-higher-rank} we deduce that the groups
        \begin{equation*}
        G_{\mathrm{II}} = \prod_{i=p+1}^{p+q} G_i \quad \text{and} \quad
        G'_{\mathrm{II}} = \prod_{j=p'+1}^{p'+q'} G'_j,
    \end{equation*}
    are isomorphic.
    We now claim that $\sigma_{\mathrm{I}}(\{ 1,\ldots, p \}) = \{ 1,\ldots, p' \}$, especially $p=p'$. If it was not the case, then some $H_i$, say $H_{i_0}$, with $i_0>p$, would be $O(u)$-equivalent to a rank one symmetric space. But then the group $\operatorname{SBE}(H_{i_0})^{O(u \vee \log)}$ would be transitive on its Gromov boundary. This cannot be, as $H_{i_0}$ has the strong pointed sphere property.
    It now follows from \eqref{item:coarse-typeI} and \cite{PalSBErankone} that the factors of 
    \begin{equation*}
        G_{\mathrm{I}} = \prod_{i=1}^{p} G_i \quad \text{and} \quad
        G'_{\mathrm{I}} = \prod_{j=1}^{p} G'_j,
    \end{equation*}
    are pairwise isomorphic, so that $G_{\mathrm{I}}$ and $G'_{\mathrm{I}}$ are isomorphic.
    Finally, by the claim above $\sigma_{\mathrm{I}}(\{ p+1,\ldots, \ell \}) = \{ p+1,\ldots, \ell \}$ so that  $H_i$ is $O(u \vee \log)$-bilipschitz equivalent to $H'_{\sigma(i)}$ for all $ i \in \{p+1,\ldots, \ell\}$. 
\end{proof}

\section{Applications to $5$-dimensional solvable Lie groups}\label{sec: contribution of product theorem}
We recall our motivation for the product theorem and the strategy of our work. Let $G$ and $H$ be two completely solvable Lie groups which we want to determine whether they are quasiisometric or not. If they were quasiisometric, we would have obtained an $O(\log)$-bilipschitz equivalence between $\rho_1(G)$ and $\rho_1(H)$. If it so happens that $\rho_1(G)$ or $\rho_1(H)$ decompose as products, and if the factors of these products admit the conditions of Theorem~\ref{thm:KKLgeneralized}, we obtain an $O(\log)$-bilipschitz equivalence between each of the respective factors. In favorable situations, we are able to rule out these factor maps by various reasons which were not applicable to the product groups. This rules out the existence of the original quasiisometry between $G$ and $H$. 

To conclude, our strategy is beneficial whenever: 
\begin{enumerate}
    \item $\rho_1(G)$ or $\rho_1(H)$ decompose as direct products, whereas $G$ and $H$ did not. 
    \item There is some obstruction for $O(\log)$-bilipschitz equivalence between the factors of $\rho_1(G)$ and $\rho_1(H)$.
\end{enumerate}

An important quasiisometry invariant is the Dehn function. 
In a subsequent paper~\cite{dehnlow}, using the work of Cornulier and Tessera~\cite{CoTesDehn}, we compute the Dehn functions of all completely solvable groups of exponential growth up to dimension 5. Based on these computations, we can summarize the contribution of our strategy to these groups. 

First, we extract the simply connected solvable Lie groups $G$ such that $\rho_1(G)$ decomposes as a product, and then divide them according to their Dehn function. The result is Table~\ref{tab:decomposable_rho1_dehn}. The names of the groups are from the classification found in \cite{Mubarakzyanov} (see also \cite{PateraZassenhaus}); the Lie groups are named $G_{d,i}^{\alpha_1,\ldots,\alpha_r}$ where $d$ is the dimension, $i$ is a positive integer, and $\alpha_1,\ldots, \alpha_r$ are parameters (we keep the parameters in the same order as in \cite{Mubarakzyanov,PateraZassenhaus}, but we sometimes use different letters, in accordance with \cite{dehnlow})).

\begin{table}[h]
    \centering
    \begin{tabular}{|c|c|c|c|c|c|}
\hline
     \multirow{2}{*}{Image by $\rho_1$} & \multicolumn{2}{c}{\textbf{Dehn function}} & & \textbf{conedim} \\
        \cline{2-4}
        & {exponential} & {quadratic} & {cubic} & \\
        \hline
    \multirow{6}{*}{}      $\mathbf{R}^2 \times G_{3,3}$ &  & $G_{5,16}^{0, \tau}$, $G_{5,17}^{\tau, 0, 1}$& & 3 \\
          $\mathbf{R}^2 \times G_{3,5}^{1/\alpha}$ & $G_{5,13}^{\alpha<1, 0, 1}$ & $G_{5,13}^{\alpha>1, 0, 1}$ & & 3 \\
    
        $\mathbf{R} \times G_{4,5}^{\gamma,1}$ \textsuperscript{($\dagger$)}  & $G_{5,19}^{1, \beta<0}$,  $G_{5,35}^{0, \beta<0}$ &
        $G_{5,19}^{1, \beta>0}$, $G_{5,35}^{0, \beta>0}$ $G_{5,27}$, $G_{5,28}^1$, $G_{5,32}^\alpha$ &  & 2 \\
       
        $\mathbf{R} \times G_{4,8}$ & $G_{5,20}^0$ & & & 2 \\
        $\mathbf{Heis} \times A_2$ &  & & $G_{5,25}^{1, 0}$ & 4 \\
        $\mathbf{R} \times G_{4,9}^1$ &  &  $G_{5,30}^1$, $G_{5,37}$ & & 2  \\
        \hline
    \end{tabular}

    \footnotesize{\textsuperscript{($\dagger$)} $\gamma$ is either equal to $\beta$, or to $1$ if there is no parameter called $\beta$.}

    \caption{Indecomposable, completely solvable groups $G$ such that $\rho_1(G)$ is decomposable, their Dehn functions, and the dimensions of their asymptotic cones.}
    \label{tab:decomposable_rho1_dehn}
\end{table}

The next step is to determine which of the factors that appear in the decompositions in the $\rho_1$ column admit the assumptions of Theorem~\ref{cor:KKLgeneralized} or at the very least Theorem \ref{th:technical-KKL}. We do this in the following Lemma.

\begin{lemma}\label{lem:check-assumptions}
Let $\alpha \in (0,1), \beta \in (0,1)$. The following hold true:
\begin{enumerate}
    \item The groups $G_{3,5}^{\alpha}, G_{4,5}^{\beta,1}$ are diagonal type Heintze groups with abelian derived subgroups.
    \item $P_1 = G_{3,3}$, $P_2 = G_{4,5}^{1,1}$ and $P_3 = G_{4,9}^1$ can be written as $P_i=A_iN_i$, $1\leqslant i \leqslant 3$, where $G_i=K_iA_iN_i$ is among the rank one simple Lie groups $G_1 = \mathrm{SO}(3,1)$,$G_2 = \mathrm{SO}(4,1)$ and $G_3 = \mathrm{SU}(2,1)$.
    \item 
    $G_{4,9}^\beta$ is a diagonal type Heintze group with the strong pointed sphere property.
\end{enumerate}
\end{lemma}

\begin{proof}
Let us check (1) with the help of the structure constants given in Table \ref{tab:structure} and using Heintze's characterisation \cite{Heintze}. For $G = G_{3,5}^\alpha$ we may write that the Lie algebra is $\mathfrak g_{5,3}^\alpha = \mathbf R^2 \rtimes_\delta \mathbf R$ where $\delta = \operatorname{ad}_{e_3}$ has strictly positive spectrum $\{1,\alpha\}$. Similarly $\mathfrak g_{4,5}^{\beta,1} = \mathbf R^2 \rtimes_\delta \mathbf R$ where $\delta = \operatorname{ad}_{e_4}$ has strictly positive spectrum $\{1,\beta\}$ (with the eigenvalue $1$ having multiplicity $2$).
As for (2), we see again from the structure constants that $\mathfrak g_{3,3} = \mathbf R^2 \rtimes_\delta \mathbf R$ where $\delta = \operatorname{ad}_{e_3} = 1$, $\mathfrak g_{4,5}^{1,1} = \mathbf R^3 \rtimes_\delta \mathbf R$ where $\delta = 1$, and $\mathfrak g_{4,9}^{1}= \mathfrak{heis}\rtimes_\delta \mathbf R$ where $\mathfrak{heis}$ is the $3$-dimensional Heisenberg Lie algebra and $\delta$ is its Carnot derivation (see e.g. \cite[\S 9.3]{PansuCCqi}).
As for (3), note that $\beta >0$ and $\beta \neq 1$ so that it is a diagonal Heintze group with derived subgroup isomorphic to $\mathbf {Heis}$, the foliation in its Gromov boundary along left cosets of $\langle e_2 \rangle$ are preserved by the self-sublinear bilipschitz equivalences \cite[Lemma 3.9]{pallier2019conf}, and so it has the strong pointed sphere property by using the exact same mechanism as in the proof of \cite[Lemma 4.1]{loglie}.
\end{proof}

\begin{table}[h]
    \centering
    \begin{tabular}{|l|l|}
\hline
$A_2$ & $[e_2,e_1]  = e_1$. \\
$\mathbf {Heis}$ & $[e_1,e_2]  = e_3$. \\
     $G_{3,3}$ & $[e_3,e_1]=e_1$, $[e_3,e_2]=e_2$. \\
     $G_{3,5}^\alpha$ & $[e_3,e_1]=e_1$, $[e_3,e_2]= \alpha e_2$, \footnotesize{$-1 \leqslant \alpha < 1$, $\alpha \neq 0$} \\
     $G_{4,5}^{\alpha, \beta}$ & $[e_4,e_1] = e_1$, $[e_4,e_2] = \alpha e_2$, $[e_4,e_3] = \beta e_3$, \footnotesize{$-1\leqslant \alpha \leqslant \beta \leqslant 1$, $\alpha\beta \neq 0$.} \\
     $G_{4,8}$ & $[e_1,e_2]=e_3$, $[e_4,e_1] = e_1$, $[e_4,e_2] = -e_2$
     \\
     $G_{4,9}^{\beta}$ & $[e_1,e_2] = e_3$, $[e_4,e_1] = e_1$, $[e_4,e_2] = \beta e_2$, $[e_4,e_3] = (1+\beta) e_3$, \footnotesize{$-1<\beta \leqslant 1$.} \\
     \hline 
    \end{tabular}

    \caption{Structure constants for the factor groups appearing in $\rho_1(G)$, $G$ from Table \ref{tab:decomposable_rho1_dehn} (whenever we do not write a bracket, it means it is $0$).}
    \label{tab:structure}
\end{table}

We conclude that if a group $G$ appearing in Table~\ref{tab:decomposable_rho1_dehn} has quadratic Dehn function then $\rho_1(G)$ admits the hypothesis of Theorem~\ref{cor:KKLgeneralized}.

\begin{corollary}[of Theorem \ref{cor:KKLgeneralized}] \label{cor: contribution of product theorem}
Denote:
\begin{itemize}
    \item $\mathcal{G}^3_{3,3}:=\{G_{5,16}^{0,\tau},G_{5,17}^{\tau,0,1}\}$
    \item$\mathcal{G}^3_{3,5}:=\{G_{5,13}^{\alpha>1,0,1}\}$
    \item $\mathcal{G}^2_{4,5}:=\{G_{5,19}^{1, \beta>0},G_{5,35}^{0, \beta>0}, G_{5,27}, G_{5,28}^1, G_{5,32}^\alpha\}$
    \item $\mathcal{G}^2_{4,9}=\{G_{5,30}^1,G_{5,37}\}$
 
\end{itemize}
Let $c \in \{2,3\}$.
    If $G\in \mathcal{G}_{i,j}^c, H\in\mathcal{G}_{k,l}^c$ are quasi-isometric, then $(i,j)=(k,l)$. Further,
    $\rho_1(G)$ and $\rho_1(H)$ are isomorphic.    
\end{corollary}

\begin{remark}
    The index $c$ in the statement serves to record the dimension of the asymptotic cone, so that if $G\in \mathcal{G}_{i,j}^c, H\in\mathcal{G}_{k,l}^{c'}$ are quasiisometric, then $c=c'$; we mention this in order to emphasize in which respect our corollary brings new information.
    If $(i,j)=(3,5)$, then one can further deduce that $G$ and $H$ are isomorphic, but at this stage, considering $\rho_0$ was actually sufficient to reach such a conclusion.
    All the groups in $\mathcal{G}_{3,3}^3$ are quasiisometric to each other (since they have the same image under $\rho_0$) so one should not hope to have the conclusion that $G$ and $H$ are isomorphic
    when $(i,j) = (3,5)$.
    On the other hand, it would be highly desirable to know if the conclusion of the corollary can be improved to a group isomorphism between $G$ and $H$ when $(i,j) \in \{(4,5), (4,9) \}$; in the latter case when $(i,j) = (4,9)$, this is related to Question \ref{ques:appn} in our Appendix, since $\mathbf R \times G_{4,9}^1$ has a Riemannian symmetric metric with a non-trivial Euclidean factor.
\end{remark}

\begin{corollary}\label{cor: g519beta}
The groups $G_{5,19}^{1,\beta}$, $\beta \in (0,+\infty)$, whose Lie algebra have structure constants 
\[ [e_1,e_2]=e_3,\, [e_5,e_1]=e_1,\, [e_5,e_3] = e_3,\, [e_5,e_4] = \beta e_4, \]
are pairwise non-quasiisometric.    
\end{corollary}

\begin{proof}
    We computed that $\rho_1(G_{5,19}^{1,\beta})$ is isomorphic to $\mathbf R \times G_{4,5}^{1,\beta}$ (see Table~\ref{tab:decomposable_rho1_dehn}) so that $\rho_1(G)$ records the parameter $\beta$ for $G\in \mathcal G_{4,5}^2$, hence by Corollary \ref{cor: contribution of product theorem}, $\beta$ is a quasiisometry invariant in this family.
\end{proof}

\begin{remark}\label{rem:end-remark}
    Using Theorem~\ref{th:technical-KKL} and part (3) in Lemma~\ref{lem:check-assumptions}, we can prove that $\mathbf R\times G_{4,9}^\beta$ is never quasiisometric to $G_{5,30}^1$ or $G_{5,37}$ unless $\beta=1$. This was not possible with Theorem~\ref{cor:KKLgeneralized} only.
\end{remark}

\begin{appendix}
    
\section{Completely solvable groups quasiisometric to symmetric spaces}
\label{app:completely-solvable-symmetric}

The theorem below follows from the combined works of many authors on the quasiisometric rigidity of symmetric spaces in the 1990s, complemented by an improvement of Kleiner-Leeb \cite{KleinerLeebRemark} and synthetized in \cite[Theorem 19.25]{CornulierQIHLC}.

\begin{theorem}
\label{th:qi-rigidity-completely-solvable}
    Let $G$ and $H$ be two completely solvable groups. 
    Assume that
    \begin{enumerate}
        \item $G$ and $H$ are quasiisometric, and
        \item $G$ or $H$ admits a symmetric left-invariant Riemannian metric with no Euclidean factor.
    \end{enumerate}
    Then $G$ and $H$ are isomorphic. 
\end{theorem}

Note that Theorem~\ref{th:qi-rigidity-completely-solvable} without assumption (2) would be \cite[Conjecture 19.113]{CornulierQIHLC}.

\begin{proof}
    Without loss of generality, we can assume that $H$ admits a left-invariant symmetric metric $g$, so that $(H,g)$ is isometric to the symmetric space $X$.
    Since $H$ is completely solvable, $X$ must be of non-compact type. $H$ acts simply transitively by isometries on $X$, so any larger connected Lie group $H'$ of isometries of $X$ will contain non-trivial point stabilizers; the latter are compact, and a completely solvable group does not have non-trivial compact subgroups. So $H$ is maximal among the completely solvable groups of isometries of $X$. 
    Now, $G$ is quasiisometric to $X$, therefore by \cite[Theorem 19.25]{CornulierQIHLC} it has a continuous, proper, cocompact action by isometries on $X$.
    The kernel of this action is a compact, therefore trivial, subgroup of $G$ so that we may consider $G$ as a subgroup of $\operatorname{Isom}(X)$; it is a closed subgroup by properness of the action. 
    Let $\widehat G$ be a maximal completely solvable subgroup of $\operatorname{Isom}(X)$ containing $G$. By the combination of \cite[Theorem 1.11]{GordonWilson} and \cite[Theorem 4.3]{GordonWilson}, $H$ and $\widehat G$ are isomorphic. It remains to show that $G=\widehat G$.
    By \cite{HigesPengAN} we know that 
    \[ \dim G = \operatorname{asdim}_{\operatorname{AN}} G =  \operatorname{asdim}_{\operatorname{AN}} H = \dim H = \dim \widehat G, \] so $G$ and $\widehat G$ have the same dimension. Since $G\subseteq \widehat G$ and $G$ and $\widehat G$ are both completely solvable, $G=\widehat G$.
\end{proof}

\begin{ques}\label{ques:appn}
In the theorem above, can one allow Euclidean factors in assumption (2)?
\end{ques}

The answer to this question is not obviously yes since in \cite[Lemma 19.29]{CornulierQIHLC}, the assumption that there is no Euclidean factors is necessary.
\end{appendix}

\bibliographystyle{alpha}
\bibliography{gabriel}

\end{document}